\documentclass[a4paper,12pt]{article}
\usepackage{amsmath, amsfonts, amssymb,xcolor}
\makeatletter
\newbox\bk@bxb
\newbox\bk@bxa
\newif\if@bkcont
\newcount\bk@lcnt

\def\breakboxskip{2pt}
\def\breakboxparindent{1.8em}

\def\breakbox{\vskip\breakboxskip\relax
\setbox\bk@bxb\vbox\bgroup
\advance\linewidth -2\fboxrule
\hsize\linewidth\@parboxrestore
\parindent\breakboxparindent\relax}

\def\bk@split{%
\@tempdimb\ht\bk@bxb 
\advance\@tempdimb\dp\bk@bxb
\setbox\bk@bxa\vsplit\bk@bxb to\z@ 
\setbox\bk@bxa\vbox{\unvbox\bk@bxa}
\setbox\@tempboxa\vbox{\copy\bk@bxa\copy\bk@bxb}
\advance\@tempdimb-\ht\@tempboxa
\advance\@tempdimb-\dp\@tempboxa}

\def\bk@addfsepht{%
\setbox\bk@bxa\vbox{\vskip\fboxsep\box\bk@bxa}}

\def\bk@addskipht{%
\setbox\bk@bxa\vbox{\vskip\@tempdimb\box\bk@bxa}}

\def\bk@addfsepdp{%
\@tempdima\dp\bk@bxa
\advance\@tempdima\fboxsep
\dp\bk@bxa\@tempdima}

\def\bk@addskipdp{%
\@tempdima\dp\bk@bxa
\advance\@tempdima\@tempdimb
\dp\bk@bxa\@tempdima}

\def\bk@line{%
\hbox to \linewidth{%
\hskip-2\fboxsep\vrule \@width\fboxrule\hskip.5\fboxsep\vrule \@width\fboxrule\hskip1.5\fboxsep
\box\bk@bxa\hfil
}}%

\def\endbreakbox{\egroup
\ifhmode\par\fi{\noindent\bk@lcnt\@ne
\@bkconttrue\baselineskip\z@\lineskiplimit\z@
\lineskip\z@\vfuzz\maxdimen
\bk@split\bk@addfsepht\bk@addskipdp
\ifvoid\bk@bxb 
\def\bk@fstln{\bk@addfsepdp
\hskip-\parindent\vbox{\llap{\raisebox{-2ex}{\rule{1.5\fboxsep}{\fboxrule}\hskip.5\fboxsep}}\bk@line\llap{\rule{1.5\fboxsep}{\fboxrule}\hskip.5\fboxsep}}}

\else 
\def\bk@fstln{\vbox{\llap{\raisebox{-2ex}{\rule{1.5\fboxsep}{\fboxrule}\hskip.5\fboxsep}}\bk@line}\hfil%
\advance\bk@lcnt\@ne
\loop
\bk@split\bk@addskipdp\leavevmode
\ifvoid\bk@bxb 
\@bkcontfalse\bk@addfsepdp
\vtop{\bk@line\noindent\hskip-2\fboxsep{\rule{1.5\fboxsep}{\fboxrule}}}%

\else 
\bk@line
\fi
\hfil\advance\bk@lcnt\@ne
\if@bkcont\repeat}%
\fi
\leavevmode\bk@fstln\par}\vskip\breakboxskip\relax}

\def\estliea{\raisebox{.5ex}{\;\rule{2.5ex}{.2ex}\;}}

\newcommand{\bizlie}[1]{\mathop{\,\raisebox{-.5ex}{$\widehat{\raisebox{.9ex}{\rule{2.5ex}{.07ex}}}_{#1}$}\,}}

\newcommand{\estliemod}[1]{\;\raisebox{.5ex}{\rule{2.5ex}{.2ex}}_{#1}\;}

\def\smp{\smallskip\par}
\def\un{{\bf 1}}
\def\zero{\{0\}}
\def\pf{\noindent{\bf Proof~:}\ }
\def\findemo{~\leaders\hbox to 1em{\hss\  \hss}\hfill~\raisebox{.5ex}{\framebox[1ex]{}}\smp}
\def\spn{\bigskip\par\noindent}
\def\mpn{\medskip\par\noindent}
\def\smpn{\smallskip\par\noindent}
\def\normal{\mathop{\trianglelefteq}}

\def\smpn{\smallskip\par\noindent}

\def\mpoint{\;\;.}
\def\mvirg{\;\;,}

\def\Res{{\rm Res}}
\def\Ind{{\rm Ind}}
\def\Inf{{\rm Inf}}
\def\Def{{\rm Def}}
\def\Iso{{\rm Iso}}

\def\Indinf{{\rm Indinf}}
\def\Defres{{\rm Defres}}

\def\Hom{{\rm Hom}}
\def\End{{\rm End}}

\def\Inf{{\rm Inf}}
\def\Im{{\rm Im}}

\def\Aut{{\rm Aut}}
\def\Ker{{\rm Ker}}

\def\Id{{\rm Id}}

\def\op{^{op}}
\def\dsp{\displaystyle}
\def\Z{\mathbb{Z}}
\def\N{\mathbb{N}}
\def\F{\mathbb{F}}
\def\Q{\mathbb{Q}}

\def\C{\mathbb{C}}

\newcommand{\dirsum}[1]{\mathop{\oplus}_{#1}\limits}
\newcommand{\romain}[1]{\uppercase\expandafter{\romannumeral #1}}

\newcommand{\sou}[1]{\underline{#1}}
\newcommand{\sur}[1]{\,\overline{\! #1}}

\newcommand{\oplusb}[2]{\mathop{\bigoplus}_{{\scriptstyle #1}\atop{\scriptstyle #2}}}

\newenvironment{enonce}[1]{\pagebreak[2]\refstepcounter{subsection}\refstepcounter{prop}\smpn{{\bf \thesection.\arabic{prop}.\ \ #1~:}}\begin{it} }{\end{it}\smp}
\newenvironment{enonce*}[1]{\pagebreak[2]\smpn{#1~:}\begin{it} }{\end{it}\smp}
\newcommand{\result}[1]{\begin{enonce}{#1}}
\def\fresult{\end{enonce}}
\newcommand{\npar}{\smallskip\par\noindent\pagebreak[2]\refstepcounter{subsection}\refstepcounter{prop}{\bf \thesection.\arabic{prop}.\ \ }}
\newcommand{\masubsect}[1]{\medskip\par\noindent\pagebreak[3]\refstepcounter{subsection}\refstepcounter{prop}{\bf \thesection.\arabic{prop}.\ \ #1.\ }}

\newenvironment{mth}[1]{\begin{breakbox}\begin{enonce}{#1}}{\end{enonce}\end{breakbox}}
\newenvironment{mth*}[1]{\begin{breakbox}\begin{enonce*}{#1}}{\end{enonce*}\end{breakbox}}
\newenvironment{rem}[1]{\refstepcounter{subsection}\refstepcounter{prop} \mpn{{\bf \thesection.\arabic{prop}.}\ \ \bf#1\ :}}{\smp}

\def\dom{\backslash}
\makeatletter
\renewenvironment{enumerate}{\ifnum \@enumdepth >3 \@toodeep\else
      \advance\@enumdepth \@ne
      \edef\@enumctr{enum\romannumeral\the\@enumdepth}\list
      {\csname label\@enumctr\endcsname}{\setlength{\topsep}{1ex}\setlength{\itemsep}{0pt}\usecounter
        {\@enumctr}\def\makelabel##1{\hss\llap{##1}}}\fi}{\endlist}
\renewenvironment{itemize}{\ifnum \@itemdepth >3 \@toodeep\else \advance\@itemdepth \@ne
\edef\@itemitem{labelitem\romannumeral\the\@itemdepth}%
\list{\csname\@itemitem\endcsname}{\setlength{\topsep}{1ex}\setlength{\itemsep}{0pt}\def\makelabel##1{\hss\llap{##1}}}\fi}
{\endlist}
\def\@sect#1#2#3#4#5#6[#7]#8{\ifnum #2>\c@secnumdepth
    \let\@svsec\@empty\else
    \refstepcounter{#1}\edef\@svsec{\csname the#1\endcsname .\hskip .5em}\fi
    \@tempskipa #5\relax
     \ifdim \@tempskipa>\z@
       \begingroup #6\relax
         \@hangfrom{\hskip #3\relax\@svsec}{\interlinepenalty \@M #8\par}%
       \endgroup
      \csname #1mark\endcsname{#7}\addcontentsline
        {toc}{#1}{\ifnum #2>\c@secnumdepth \else
                     \protect\numberline{\csname the#1\endcsname}\fi
                   #7}\else
       \def\@svsechd{#6\hskip #3\relax  
                  \@svsec #8\csname #1mark\endcsname
                     {#7}\addcontentsline
                          {toc}{#1}{\ifnum #2>\c@secnumdepth \else
                            \protect\numberline{\csname the#1\endcsname}\fi
                      #7}}\fi
    \@xsect{#5}}
\def\section{\@startsection {section}{1}{\z@}{-3.5ex plus-1ex minus
    -.2ex}{2.3ex plus.2ex}{\reset@font\Large\bf}}  

\makeatother
\renewenvironment{equation}{\refstepcounter{subsection}\refstepcounter{prop}$$}{\leqno{\bf (\theprop)}$$}

\def\mar[#1]{\ar@{-}[#1]|-{\object@{<}}}
\def\marb[#1]{\ar@{-}[#1]|{\object+{  }}}

\usepackage[all]{xy}
\def\Inn{\rm Inn}

\def\isom{\stackrel{{_{\cong}}}{\to}}
\renewcommand{\oplusb}[2]{\mathop{\oplus}_{{\scriptstyle #1}\atop{\scriptstyle #2}}}
\def\rDelta{\overrightarrow{\Delta}}
\def\lDelta{\overleftarrow{\Delta}}

\newcommand{\monast}[2]{{{\raisebox{-.5ex}{$\scriptstyle #2$}}\atop{\rule{0ex}{4ex}\hbox{\Huge\rm *}\atop{\smash{\raisebox{1.5ex}{$\scriptstyle#1$}}}}}\!}
\def\mmp{\medskip\par}
\begin{document}
\centerline{\Large\bf The Roquette category of finite $p$-groups}\vspace{.5cm}\par
\centerline{\bf Serge Bouc
}\vspace{1cm}\par
{\footnotesize {\bf Abstract :} Let $p$ be a prime number. This paper introduces the {\em Roquette category} $\mathcal{R}_p$ of finite $p$-groups, which is an additive tensor category containing all finite $p$-groups among its objects. In $\mathcal{R}_p$, every finite $p$-group $P$ admits a canonical direct summand $\partial P$, called {\em the edge} of $P$. Moreover $P$ splits uniquely as a direct sum of edges of {\em Roquette $p$-groups}, and the tensor structure of $\mathcal{R}_p$ can be described in terms of such edges.}\par  
{\footnotesize The main motivation for considering this category is that the additive functors from~$\mathcal{R}_p$ to abelian groups are exactly the {\em rational $p$-biset functors}. This yields in particular very efficient ways of computing such functors on arbitrary $p$-groups~: this applies to the representation functors $R_K$, where $K$ is any field of characteristic 0, but also to the functor of units of Burnside rings, or to the torsion part of the Dade group.
}\vspace{2ex}\par
{\footnotesize {\bf AMS Subject classification :} 18B99, 19A22, 20C99, 20J15.\vspace{2ex}}\par
\section{Introduction}
Let $p$ be a prime number. This article introduces {\em the Roquette category} $\mathcal{R}_p$ of finite $p$-groups, which is an additive tensor category with the following properties~:
\begin{itemize}
\item Every finite $p$-group can be viewed as an object of $\mathcal{R}_p$. The tensor product of two finite $p$-groups $P$ and $Q$ in $\mathcal{R}_p$ is the direct product $P\times Q$.
\item In $\mathcal{R}_p$, any finite $p$-group has a direct summand $\partial P$, called {\em the edge} of~$P$, such that
$$P\cong\dirsum{N\normal P}\partial (P/N)\mpoint$$
Moreover, if the center of $P$ is not cyclic, then $\partial P=0$.
\item In $\mathcal{R}_p$, every finite $p$-group $P$ decomposes as a direct sum
$$P\cong\dirsum{R\in \mathcal{S}}\partial R\mvirg$$
where $\mathcal{S}$ is a finite sequence of {\em Roquette groups}, i.e. of $p$-groups of normal $p$-rank 1, and such a decomposition is essentially unique. Given the group~$P$, such a decomposition can be obtained explicitly from the knowledge of a {\em genetic basis} of $P$ (Theorem~\ref{Roquette category} and Proposition~\ref{isomorphism}).
\item The tensor product $\partial P\times\partial Q$ of the edges of two Roquette $p$-groups $P$ and $Q$ is isomorphic to a direct sum of a certain number $\nu_{P,Q}$ of copies of the edge $\partial (P\diamond Q)$ of another Roquette group (where both $\nu_{P,Q}$ and $P\diamond Q$ are known explicitly - see Theorem~\ref{structure diamond} and Corollary~\ref{edge product}).
\item The additive functors from $\mathcal{R}_p$ to the category of abelian groups are exactly the {\em rational $p$-biset functors} introduced in~\cite{bisetsections}. 
\end{itemize}
The latter is the main motivation for considering this category~: any structural result on $\mathcal{R}_p$ will provide for free some information on such rational functors for $p$-groups, e.g. the representation functors $R_K$, where $K$ is a field of characteristic 0 (see \cite{doublact}, \cite{fonctrq}, and L. Barker's article~\cite{rhetoric}), the functor of units of Burnside rings~(\cite{burnsideunits}), or the torsion part of the Dade group~(\cite{dadegroup}). \par
In particular, the above results on $\mathcal{R}_p$ yield isomorphisms describing the structure of some $p$-groups as objects of this category, and this is enough to compute the evaluations of rational $p$-biset functors. For example 
$$(D_8)^n\cong \un\oplus (5^n-1)\cdot\partial C_2$$
in $\mathcal{R}_2$ (Equation~\ref{D2n}). More generally, Proposition~\ref{D2nk} gives a formula for $(D_{2^m})^n$. A straightforward consequence, applying the functor $R_\Q$, is the following
\begin{mth*}{{\bf Example}} For any $n\in N$, the group $(D_8)^n$ has $5^n$ conjugacy classes of cyclic subgroups.
\end{mth*}
Another important by-product of the above result giving the tensor structure of $\mathcal{R}_p$ is the explicit description of a genetic basis of a direct product $P\times Q$, in terms of a genetic basis of $P$ and a genetic basis of $Q$ (Theorem~\ref{genetic basis product}). This allows in particular for a quick computation of the torsion part of the Dade group of some $p$-groups, e.g. (Theorem~\ref{product of dihedral groups}, Assertion~1 and Assertion~3)~:
\begin{mth*}{{\bf Theorem}} \begin{itemize}
\item Let $P$ be an arbitrary finite direct product of groups of order 2 and dihedral 2-groups. Then the Dade group of any factor group of $P$ is torsion free.
\item Let $n$ be a positive integer. For any integer $m\geq 4$, let $P=SD_{2^m}$ be a semidihedral group of order $2^m$, and let $P^{*n}$ denote the central product of $n$-copies of $P$. Then the torsion part of the Dade group of $P^{*n}$ is isomorphic to $(\Z/2\Z)^{2^{(n-1)(m-3)}}$.
\end{itemize}
\end{mth*}
This also yields similar results on groups of units of Burnside rings of these groups (Remark~\ref{units burnside}), or on representations of central products of $p$-groups, as in Examples~\ref{edges} and~\ref{extraspecial}~:
\begin{mth*}{\bf Example} Let $p$ be a prime, let $X$ be an extraspecial $p$-group, and let~$Q$ be a non-trivial $p$-group. Let $K$ be a field of characteristic 0. Then $Q$ has the same number (possibly 0) of isomorphism classes of faithful irreducible representations over $K$ as any central product $X*Q$.
\end{mth*}
Another possibly interesting phenomenon is that some non-isomorphic $p$-groups may become isomorphic in the category $\mathcal{R}_p$. This means that some non-isomorphic $p$-groups cannot be distinguished using only rational $p$-biset functors. When $p=2$, there are even examples where this occurs for groups of different orders (Example~\ref{different orders}). When $p>2$, saying that the $p$-groups $P$ and $Q$ are isomorphic in $\mathcal{R}_p$ is equivalent to saying that the group algebras $\Q P$ and $\Q Q$ have isomorphic centers (Proposition~\ref{isomorphism Roquette}).\par
The category $\mathcal{R}_p$ is built as follows~: consider first the category $\mathcal{R}_p^\sharp$, which is the quotient category of the biset category of finite $p$-groups (in which objects are finite $p$-groups and morphisms are virtual bisets) obtained by killing a specific element $\delta$ in the Burnside group of the Sylow $p$-subgroup of $PGL(3,\F_p)$. Then take idempotent completion, and additive completion of the resulting category.\par
In particular, this construction relies on bisets, and related functors. Consequently, the paper is organized as follows~: Section~\ref{rational p-biset functors} is a (not so) quick summary of the background on biset functors, Roquette groups, genetic bases of $p$-groups, and rational $p$-biset functors. The category $\mathcal{R}_p$ is introduced in Section~\ref{section Roquette category}, and in Section~\ref{the tensor structure}, its tensor structure is described. Finally Section~\ref{examples and applications} gives some examples and applications.
\spn
{\bf Acknowledgments :} 
Even though the idea of considering the category $\mathcal{R}_p$ was implicit in~\cite{rationnel}, the present paper wouldn't exist without the illuminating conversations I had with Paul Balmer in november 2010, during which I gradually understood I was a kind of {\em Monsieur Jourdain} of tensor categories\ldots ~I~wish to thank him for this revelation.
\section{Rational $p$-biset functors}\label{rational p-biset functors}
\masubsect{Biset functors} The {\em biset category} $\mathcal{C}$ of finite groups is defined as follows~:
\begin{itemize}
\item The objects of $\mathcal{C}$ are the finite groups.
\item Let $G$ and $H$ be finite groups. Then
$$\Hom_{\mathcal{C}}(G,H)=B(H,G)\mvirg$$
where $B(H,G)$ denotes the Grothendieck group of the category of finite $(H,G)$-bisets, i.e. the Burnside group of the group $H\times G\op$.
\item Let $G$, $H$, and $K$ be finite groups. The composition of morphisms 
$$B(K,H)\times B(H,G)\to B(K,G)$$
in the category $\mathcal{C}$ is the linear extension of the product induced by the product of bisets $(V,U)\mapsto V\times_HU$, where $V$ is a $(K,H)$-biset, and $U$ is an $(H,G)$-biset.
\item The identity morphism of the finite group $G$ is the image in $B(G,G)$ of the set $G$, endowed with its $(G,G)$-biset structure given by left and right multiplication.
\end{itemize}
\begin{mth}{Definition} A {\em biset functor} is an additive functor from $\mathcal{C}$ to the category of abelian groups. A {\em morphism} of biset functors is a natural transformation of functors.
\end{mth}
Morphisms of biset functors can be composed, and the resulting category of biset functors is denoted by $\mathcal{F}$. It is an abelian category.
\begin{rem}{Examples} \label{basic bisets}\begin{enumerate}
\item The correspondence $B$ sending a finite group $G$ to its Burnside group $B(G)$ is a biset functor, called the {\em Burnside functor}~: indeed $B(G)=\Hom_{\mathcal{C}}(\un,G)$, so $B$ is in fact the Yoneda functor $\Hom_{\mathcal{C}}(\un,{-})$.
\item The formalism of bisets gives a single framework for the usual operations of induction, restriction, inflation, deflation, and transport by isomorphism via the following correspondences~: 
\begin{itemize} 
\item If $H$ is a subgroup of $G$, then let $\Ind_H^G\in B(G,H)$ denote the set~$G$, with left action of $G$ and right action of $H$ by multiplication.
\item If $H$ is a subgroup of $G$, then let $\Res_H^G\in B(H,G)$ denote the set~$G$, with left action of $H$ and right action of $G$ by multiplication.
\item If $N\normal G$, and $H=G/N$, then let $\Inf_H^G\in B(G,H)$ denote the set $H$, with left action of~$G$ by projection and multiplication, and right action of $H$ by multiplication.
\item If $N\normal G$, and $H=G/N$, then let $\Def_H^G\in B(H,G)$ denote the
 set~$H$, with left action of $H$ by multiplication, and right action of~$G$ by  projection and multiplication.
\item If $\varphi: G\to H$ is a group isomorphism, then let $\Iso_G^H=\Iso_G^H(\varphi)\in  B(H,G)$ denote the set $H$, with left action of $H$ by multiplication, and right action of $G$ by taking image by $\varphi$, and then multiplying in $H$.
\item When $H$ is a subgroup of $G$, let $\Defres_{N_G(H)/H}^G\in B\big(N_G(H)/H,G\big)$ denote the set $H\dom G$, viewed as a $\big(N_G(H)/H,G\big)$-biset. It is equal to the composition $\Def_{N_G(H)/H}^G\circ \Res_{N_G(H)}^G$.
\item When $H$ is a subgroup of $G$, let $\Indinf_{N_G(H)/H}^G\in B\big(G,N_G(H)/H\big)$ denote the set $G/H$, viewed as a $\big(G,N_G(H)/H\big)$-biset. It is equal to the composition $\Ind_{N_G(H)}^G\circ \Inf_{N_G(H)/H}^{N_G(H)}$.
\end{itemize}
\end{enumerate}
\end{rem}
\vspace{2ex}
\masubsect{$p$-biset functors} From now on, the symbol $p$ will denote a prime number. 
\begin{mth}{Definition and Notation} 
\vspace{2ex}
\begin{itemize}
\item The {\em biset category} $\mathcal{C}_p$ of finite $p$-groups is the full subcategory of $\mathcal{C}$ consisting of finite $p$-groups.
\item A $p$-biset functor is an additive functor from $\mathcal{C}_p$ to the category of abelian groups. A morphism of $p$-biset functors is a natural transformation of functors.
\item $p$-biset functors form an abelian category $\mathcal{F}_p$.
\end{itemize}
\end{mth}
\vspace{2ex}
\pagebreak[3]
\masubsect{Roquette $p$-groups}
\begin{mth}{Definition} A finite group $G$ is call a {\em Roquette group} if it has normal rank 1, i.e. if all the normal abelian subgroups of $G$ are cyclic.
\end{mth}
The Roquette $p$-groups have been first classified by\ldots Roquette (\cite{roquette}, see also~\cite{gorenstein})~: these are the cyclic groups, if $p>2$. The Roquette 2-groups are the cyclic groups, the generalized quaternion groups, the dihedral and semidihedral groups of order at least 16. \par
More generally, the $p$-hyperelementary Roquette groups have been classified by Hambleton, Taylor and Williams (Theorem~3.A.6 of \cite{htw}). \par
The following schematic diagram represents the lattice of subgroups of the dihedral group $D_{16}$, the quaternion group $Q_{16}$, and the semi-dihedral group $SD_{16}$ (an horizontal dotted link between two vertices means that the corresponding subgroups are conjugate)~:
\def\xar[#1]{\ar@{-}[#1]}
\newcommand\boitext[1]{\makebox[0pt][l]{$\scriptstyle #1$}}
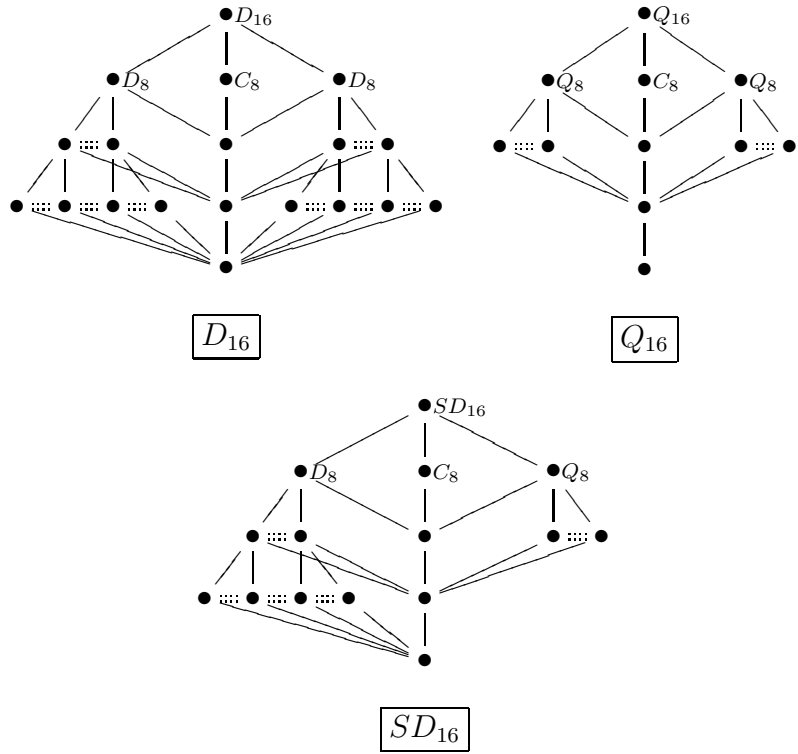
$$\xymatrix@C=6pt@R=12pt{
&&&&\bullet\boitext{D_{16}}&&&&\\
&&\bullet\xar[urr]\boitext{D_8}&&\bullet\xar[u]\boitext{C_8}&&\bullet\xar[ull]\boitext{D_8}&\\
&\bullet\xar[ur]\ar@{::}[r]&\bullet\xar[u] &&\bullet\xar[ull]\xar[u]\xar[urr]&&\bullet\xar[u]&\bullet\xar[ul]\ar@{::}[l]\\
\bullet\xar[ur]\ar@{::}[r]&\bullet\xar[u]\ar@{::}[r]&\bullet\xar[u]\ar@{::}[r]&\bullet\xar[ul]&\bullet\xar[ull]\xar[u]\xar[urr]\xar[ulll]\xar[urrr]&\bullet\xar[ur]&\bullet\xar[u]\ar@{::}[l]&\bullet\xar[u]\ar@{::}[l]&\bullet\xar[ul]\ar@{::}[l]\\
&&&&\bullet\xar[ullll]\xar[ulll]\xar[ull]\xar[ul]\xar[u]\xar[ur]\xar[urr]\xar[urrr]\xar[urrrr]\\
&&&&*+[F]{D_{16}}\\
}
\xymatrix@C=6pt@R=12pt{
&&&&\bullet\boitext{Q_{16}}&&&&\\
&&\bullet\xar[urr]\boitext{Q_8}&&\bullet\xar[u]\boitext{C_8}&&\bullet\xar[ull]\boitext{Q_8}&\\
&\bullet\xar[ur]\ar@{::}[r]&\bullet\xar[u] &&\bullet\xar[ull]\xar[u]\xar[urr]&&\bullet\xar[u]&\bullet\xar[ul]\ar@{::}[l]\\
&&&&\bullet\xar[ull]\xar[u]\xar[urr]\xar[ulll]\xar[urrr]&&&&\\
&&&&\bullet\xar[u]\\
&&&&*+[F]{Q_{16}}\\
}$$
$$
\xymatrix@C=6pt@R=12pt{
&&&&\bullet\boitext{SD_{16}}&&&&\\
&&\bullet\xar[urr]\boitext{D_8}&&\bullet\xar[u]\boitext{C_8}&\ \ &\bullet\xar[ull]\boitext{Q_8}&\\
&\bullet\xar[ur]\ar@{::}[r]&\bullet\xar[u] &&\bullet\xar[ull]\xar[u]\xar[urr]&&\bullet\xar[u]&\bullet\xar[ul]\ar@{::}[l]\\
\bullet\xar[ur]\ar@{::}[r]&\bullet\xar[u]\ar@{::}[r]&\bullet\xar[u]\ar@{::}[r]&\bullet\xar[ul]&\bullet\xar[ull]\xar[u]\xar[urr]\xar[ulll]\xar[urrr]&&&&\\
&&&&\bullet\xar[ullll]\xar[ulll]\xar[ull]\xar[ul]\xar[u]\\
&&&&*+[F]{SD_{16}}\\
}$$
These diagrams give a good idea of the general case. 
\begin{mth}{Definition} \label{axial}Let $G$ be a finite group, of exponent $e$. An {\em axis} of $G$ is a cyclic subgroup of order $e$ in $G$. An {\em axial} subgroup of $G$ is a subgroup of an axis of $G$.
\end{mth}
With these definitions, let us recall without proof the following properties of Roquette $p$-groups~:
\begin{mth}{Lemma}\label{en vrac} Let $P$ be a non-trivial Roquette $p$-group, of exponent $e_P$.
\begin{enumerate}
\item The center of $P$ is cyclic, hence $P$ admits a unique central subgroup $Z_P$ of order $p$.
\item There exists a non-trivial subgroup $Q$ of $P$ such that $Q\cap Z(P)=\un$ if and only if $p=2$, and $P$ is dihedral or semidihedral. In this case moreover $|Q|=2$, and $N_P(Q)=QZ_P$.
\item If $P$ is not cyclic, then $p=2$ and $e_P=|P|/2$.
\item There is a unique axis in $P$, except in the case $P\cong Q_8$, where there are three of them. Any axis of $P$ is normal in  $P$.
\item If $R$ is a non-trivial axial subgroup of $P$, then $R\geq Z_P$, and $R\normal P$. If moreover $|R|\geq p^2$, then $C_P(R)$ is the only axis of $P$ containing $R$. 
\end{enumerate}
\end{mth}
Let us also recall the following~:
\begin{mth}{Lemma} Let $P$ be a finite Roquette $p$-group. Then there is a unique simple faithful $\Q P$-module $\Phi_P$, up to isomorphism.
\end{mth}
\pf See \cite{bisetfunctors} Proposition 9.3.5.\findemo
\begin{rem}{Example} \label{endo Phi}Let $P$ be a cyclic group of order $p^m$, and suppose first that $m>1$. The algebra $\Q P$ is isomorphic to the algebra $A=\Q[X]/(X^{p^m}-1)$. As $$X^{p^m}-1=(X^{p^{m-1}}-1)\Psi_{p^m}(X)\mvirg$$
where $\Psi_{p^m}$ denotes the $p^m$-th cyclotomic polynomial, it follows that there is a split exact sequence of $A$-modules
$$0\mapsto \Q[X]/(\Psi_{p^m})\to A\to \Q[X]/(X^{p^{m-1}}-1)\to 0\mvirg$$
which can be viewed as the following sequence of $\Q P$-modules
$$0\to \Phi_P\to \Q P\to \Q (P/Z)\to 0\mvirg$$
where $Z$ is the unique subgroup of order $p$ of $P$. It follows that there is an isomorphism of $\Q$-algebras
$$\End_{\Q P}\Phi_P\cong \Q[X]/(\Psi_{p^m})\cong\Q(\zeta_{p^m})\mvirg$$
where $\zeta_{p^m}$ is a primitive $p^m$-th root of unity in $\C$. \par
Now if $m=0$, then $P=\un$, $\Phi_P=\Q$, and $\End_{\Q P}\Phi_P\cong \Q$, also.
\end{rem}
\masubsect{Expansive and genetic subgroups}
\begin{mth}{Definition} A subgroup $H$ of a group $G$ is called {\em expansive} if for any $g\in G$ such that $H^g\neq H$, the group $\big(H^g\cap N_G(H)\big)H/H$ contains a non trivial normal subgroup of $N_G(H)/H$, i.e. if
$$g\in G-N_G(H)\Longrightarrow \bigcap_{n\in N_G(H)}\big(H^{gn}\cap N_G(H)\big)H>H\mpoint$$
\end{mth}
\begin{rem}{Example} \label{normalisateur normal}If $H\normal G$, then $H$ is expansive in $G$. More generally, if $N_G(H)\normal G$, then $H$ is expansive in $G$~: indeed $N_G(H^g)=N_G(H)$, for any $g\in G$. Hence for $g\in G-N_G(H)$
$$\bigcap_{n\in N_G(H)}\big(H^{gn}\cap N_G(H)\big)H=\big(H^g\cap N_G(H)\big)H=H^g\cdot H>H\mpoint$$
\end{rem}
\begin{mth}{Notation} When $H$ is a subgroup of the group $G$, denote by $Z_G(H)$ the subgroup of $N_G(H)$, containing $H$, defined by
$$Z_G(H)/H=Z\big(N_G(H)/H\big)\mpoint$$
\end{mth}
The following is an easy consequence of well known properties of $p$-groups~:
\begin{mth}{Lemma} {\rm [\cite{bisetfunctors} Lemma 9.5.2]} Let $Q$ be a subgroup of a finite $p$-group~$P$. Then $Q$ is expansive in $P$ if and only if
$$\forall g\in P,\;\;Q^g\cap Z_P(Q)\leq Q\Longrightarrow Q^g=Q\mpoint$$
\end{mth}
\begin{mth}{Example} \label{diagonal}Let $P$ be a $p$-group, and let 
$$\Delta(P)=\big\{(x,x)\in (P\times P)\mid x\in P\big\}$$
denote the diagonal subgroup of $(P\times P)$. Then $N_{P\times P}\big(\Delta(P)\big)/\Delta(P)\cong Z(P)$, and $\Delta(P)$ is expansive in $(P\times P)$ if and only if
$$\forall x\in P,\;\;[P,x]\cap Z(P)=\{1\}\implies x\in Z(P)\mvirg$$
where $[P,x]$ is the {\em set} of commutators $[y,x]=y^{-1}x^{-1}yx=y^{-1}y^x$, for $y\in P$.\par
In particular, if $[P,P]\leq Z(P)$, then $\Delta(P)$ is expansive in $(P\times P)$. 
\end{mth}
\pf Indeed $N_{P\times P}\big(\Delta(P)\big)$ consists of pairs $(a,b)\in (P\times P)$ such that $ab^{-1}\in Z(P)$. This shows that $N_{P\times P}\big(\Delta(P)\big)/\Delta(P)\cong Z(P)$, and that 
$$Z_{P\times P}\big(\Delta(P)\big)=N_{P\times P}\big(\Delta(P)\big)=\big(\un\times Z(P)\big)\Delta(P)\mpoint$$
Now $\Delta(P)^{(u,v)}=\Delta(P)^{(1,x)}$, for any $(u,v)\in (P\times P)$, where $x=u^{-1}v$, and
$$\Delta(P)^{(1,x)}\cap Z_{P\times P}\big(\Delta(P)\big)=\big\{(t,t^x)\mid t\in P,\;t^{-1}t^x\in Z(P)\big\}\mpoint$$
Hence $\Delta(P)^{(1,x)}\cap Z_{P\times P}\big(\Delta(P)\big)\leq \Delta(P)$ if and only if for any $t\in P$, the assumption $t^{-1}t^x\in Z(P)$ implies $t=t^x$, i.e. $[t,x]=1$, in other words if $[P,x]\cap Z(P)=\{1\}$. Hence $\Delta(P)$ is expansive in $(P\times P)$ if and only if for any $x\in P$, the assumption $[P,x]\cap Z(P)=\{1\}$ implies $(1,x)\in N_{P\times P}\big(\Delta(P)\big)$, i.e. $x\in Z(P)$, as claimed. The last assertion follows trivially.\findemo
\begin{mth}{Definition} Let $Q$ be a subgroup of the finite $p$-group $P$. Then $Q$ is called a {\em genetic} subgroup of $P$ if $Q$ is expansive in $P$ and if the group $N_P(Q)/Q$ is a Roquette group.
\end{mth}
\begin{mth}{Definition} Define a relation $\bizlie{P}$ on the set of subgroups of the finite $p$-group $P$ by
$$Q\bizlie{P}R\;\;\Leftrightarrow \;\;\exists g\in P,\;\;Q^g\cap Z_P(R)\leq R\;\;\hbox{and}\;\;{^gR}\cap Z_P(Q)\leq Q\mpoint$$\vspace{-4ex}
\end{mth}
\begin{mth}{Lemma} \label{normalisateur normal bizlie}Let $P$ be a finite $p$-group. If $Q$ and $R$ are subgroups of $P$ such that $N_P(Q)=N_P(R)\normal P$, then
$$Q\bizlie{P}R\;\;\Leftrightarrow\;\; Q=_PR\mpoint$$
\end{mth}
\pf Indeed, since $P$ is a $p$-group, saying that $Q^g\cap Z_P(R)\leq R$ is equivalent to saying that the subgroup $\big(Q^g\cap N_P(R)\big)R/R$ of $N_P(R)/R$ contains no non-trivial normal subgroup of $N_P(R)/R$, i.e. that
$$\bigcap_{n\in N_P(R)}\big(Q^{gn}\cap N_P(R)\big)R=R\mpoint$$
But if $N_P(Q)\normal P$, then $N_P(Q)=N_P(Q^g)=N_P(R)$, for any $g\in G$. Hence
$$\bigcap_{n\in N_G(R)}\big(Q^{gn}\cap N_P(R)\big)R=Q^g\cdot R\mpoint$$
This is equal to $R$ if and only if $Q^g\leq R$. Similarly ${^gR}\cap Z_P(Q)\leq Q$ if and only if $^gR\leq Q$. Hence $Q^g=R$.\findemo
\begin{mth}{Definition} Let $G$ be a group.
\begin{enumerate}
\item A {\em section} of $G$ is a pair $(T,S)$ of subgroups of $G$ such that $S\normal T$. The quotient $T/S$ is called the corresponding {\em subquotient} of $G$.
\item Two sections $(T,S)$ and $(Y,X)$ of $G$ are said to be {\em linked} (notation $(T,S)\estliea (Y,X)$) if 
$$S(T\cap Y)=T,\;\;X(T\cap Y)=Y,\;\;T\cap X=S\cap Y\mpoint$$
They are said to be {\em linked modulo $G$} (notation $(T,S)\estliemod{G}(Y,X)$) if there exists $g\in G$ such that $(T,S)\estliea ({^gY},{^gX})$.
\end{enumerate}
\end{mth}
Observe in particular that if $(T,S)\estliemod{G}(Y,X)$, then the corresponding subquotients $T/S$ and $Y/X$ are isomorphic.

\begin{mth}{Theorem} \label{genetic recall}Let $P$ be a finite $p$-group. 
\begin{enumerate}
\item If $S$ is a genetic subgroup of $P$, then the module
$$V(S)=\Ind_{N_P(S)}^P\Inf_{N_P(S)/S}^{N_P(S)}\Phi_{N_P(S)/S}$$
is a simple $\Q P$-module. Moreover, the functor $\Ind_{N_P(S)}^P\Inf_{N_P(S)/S}^{N_P(S)}$ induces an isomorphism of $\Q$-algebras
$$\End_{Q P}V(S)\cong\End_{\Q N_P(S)/S}\Phi_{N_P(S)/S}\mpoint$$
\item If $V$ is a simple $\Q P$-module, then there exists a genetic subgroup $S$ of $P$ such that $V\cong V(S)$.
\item If $S$ and $T$ are genetic subgroups of $P$, then
$$V(S)\cong V(T) \;\;\Leftrightarrow \;\;S\bizlie{P}T\;\;\Leftrightarrow\;\; \big(N_P(S),S\big)\estliemod{P}\big(N_P(T),T\big)\mpoint$$
In particular, if $S\bizlie{P}T$, then $N_P(S)/S\cong N_P(T)/T$. Moreover, the relation $\bizlie{P}$ is an equivalence relation on the set of genetic subgroups of $P$, and the corresponding set of equivalence classes is in one to one correspondence with the set of isomorphism classes of simple $\Q P$-modules.
\end{enumerate}
\end{mth}
\pf See \cite{bisetfunctors}, Theorem 9.6.1.\findemo
\begin{mth}{Definition} Let $P$ be a finite $p$-group. A {\em genetic basis} of $P$ is a set of representatives of equivalence classes of genetic subgroups of $P$ for the relation $\bizlie{P}$.
\end{mth}
\masubsect{Faithful elements} (\cite{bisetfunctors} Sections 6.2 and 6.3) Let $G$ be a finite group. If $N$ is a normal subgroup of $G$, recall from Examples~\ref{basic bisets} that $\Inf_{G/N}^G$ denotes the set~$G/N$, viewed as a $(G,G/N)$-biset for the actions given by (projection to the factor group and) multiplication in $G/N$. Similarly $\Def_{G/N}^G$ denotes the same set~$G/N$, considered as a $(G/N,G)$-biset.\par
There is an isomorphism of $(G/N,G/N)$-bisets
\begin{equation}\label{definf}
\Id_{G/N}\cong \Def_{G/N}^G\circ \Inf_{G/N}^G\mpoint
\end{equation}
More generally, if $M$ and $N$ are normal subgroups of $G$, there is a isomorphism of $(G/M,G/N)$-bisets
\begin{equation}\label{infdef}
\Def_{G/M}^G\circ\Inf_{G/N}^G=\Inf_{G/MN}^{G/M}\circ\Def_{G/MN}^{G/N}\mpoint
\end{equation}
It follows that if $j_N^G$ is defined by
$$j_N^G=\Inf_{G/N}^G\circ\Def_{G/N}^G\mvirg$$
then $j_M^G\circ j_N^G=j_{MN}^G$. In particular $j_N^G$ is an idempotent of $B(G,G)$. Moreover, by a standard orthogonalization procedure, the elements $f_N^G$ defined for~$N\normal G$ by
$$f_N^G=\sum_{N\leq M\normal G}\mu_{\normal G}(N,M)j_M^G\mvirg$$
where $\mu_{\normal G}(N,M)$ is the M\"obius function of the poset of normal subgroups of~$G$, are orthogonal idempotents of $B(G,G)$, and their sum is equal to $\Id_{G/N}$. The idempotent $f_\un^G$ is of special importance~:
\begin{mth}{Lemma} \label{funG}Let $G$ be a finite group, and $N$ be a normal subgroup of~$G$. 
\begin{enumerate}
\item $f_N^G=\Inf_{G/N}^G\circ f_\un^{G/N}\circ \Def_{G/N}^G$ in $B(G,G)$.
\item If $N\neq 1$, then $\Def_{G/N}^G\circ f_\un^G=0$ in $B(G/N,G)$, and $f_\un^G\circ\Inf_{G/N}^G=0$ in $B(G,G/N)$.
\end{enumerate}
\end{mth}
\pf Assertion~1 is Remark~6.2.9 of~\cite{bisetfunctors}, and Assertion~2 is a special case of Proposition~6.2.6.\findemo
If $F$ is a biset functor, the set $\partial F(G)$ of {\em faithful elements} of $F(G)$ is defined by
$$\partial F(G)=F(f_\un^G)F(G)\mpoint$$
It can be shown (\cite{bisetfunctors} Lemma 6.3.2) that 
$$\partial F(G)=\mathop{\bigcap}_{\un<N\normal G}\limits\Ker\; F(\Def_{G/N}^G)\mpoint$$
\begin{rem}{Example} \label{representation edge}Let $F=R_K$ be the representation functor over a field $K$ of characteristic 0. Then for a finite group $G$, the group $\partial R_K(G)$ is the direct summand of $R_K(G)$ with basis the set of (isomorphism classes of) faithful irreducible $KG$-modules.
\end{rem}
\begin{mth}{Lemma}\label{inflation} Let $G$ be a group, and $S$ be a subgroup of $G$, such that $S\cap Z(G)\neq \un$. Then $\Defres_{N_G(S)/S}^Gf_\un^G=0$ in $B\big(N_G(S)/S,G\big)$.
\end{mth}
\pf Set $N=S\cap Z(G)$, $\sur{G}=G/N$, and $\sur{S}=S/N$. Then
$$\Defres_{N_G(S)/S}^G=\Defres_{N_{\sur{G}}(\sur{S})/\sur{S}}^{\sur{G}}\Def_{G/N}^G\mpoint$$
Now $\Def_{G/N}^Gf_\un^G=0$ if $N\neq \un$, by Lemma~\ref{funG}. \findemo
\begin{mth}{Theorem} {\rm [\cite{bisetfunctors} Theorem 10.1.1]} \label{roquette plus}Let $P$ be a finite $p$-group, and $\mathcal{B}$ be a genetic basis of $P$. Then, for any $p$-biset functor $F$, the map
$$\mathcal{I}_\mathcal{B}=\dirsum{S\in\mathcal{B}}\Indinf_{N_P(S)/S}^P:\dirsum{S\in\mathcal{B}}\partial F\big(N_P(S)/S\big)\to F(P)$$
is split injective. A left inverse is the map
$$\mathcal{D}_\mathcal{B}=\dirsum{S\in\mathcal{B}}f_\un^{N_P(S)/S}\circ \Defres_{N_P(S)/S}^P:F(P)\to \dirsum{S\in\mathcal{B}}\partial F\big(N_P(S)/S\big)\mpoint$$
\end{mth}
One can show (\cite{bisetfunctors} Lemma 10.1.2) that if $\mathcal{B}$ and $\mathcal{B}'$ are genetic bases of $P$, the map $\mathcal{I}_\mathcal{B}$ is an isomorphism if and only if the map $\mathcal{I}_{\mathcal{B}'}$ is an isomorphism. This motivates the following definition~:
\begin{mth}{Definition} \label{rational functors} A $p$-biset functor $F$ is called {\em rational} if for any finite $p$-group $P$, there exists a genetic basis $\mathcal{B}$ of $P$ such that the map $\mathcal{I}_\mathcal{B}$ is an isomorphism.
\end{mth}
So $F$ is rational if and only if for any finite $p$-group $P$ and {\em any} genetic basis $\mathcal{B}$ of $P$, the map $\mathcal{I}_\mathcal{B}$ is an isomorphism.
\begin{rem}{Examples} \mpn
$\bullet$ The functor $R_\Q$ of rational representations, which sends the finite $p$-group $P$ to the group $R_\Q(P)$, is a rational $p$-biset functor. This example is of course the reason for calling {\em rational} the $p$-biset functors of Definition~\ref{rational functors}. This choice has proved rather unfortunate, since the $p$-biset functor $R_\C$ of {\em complex} representations is also a rational functor\ldots More generally, if $K$ is a field of characteristic 0, then the functor $R_K$ is a rational $p$-biset functor.\mpn
$\bullet$ The functor of units of the Burnside ring, sending a $p$-group $P$ to the group of units $B^\times(P)$ of its Burnside ring, is a rational $p$-biset functor (see~\cite{burnsideunits}).\mpn
$\bullet$ Let $k$ be a field of characteristic $p$. The correspondence sending a finite $p$-group $P$ to the torsion part $D_k^t(P)$ of the Dade group of $P$ over $k$ is not a biset functor in general, because of phenomenons of {\em Galois twists}, but still the maps $\mathcal{I}_{\mathcal{B}}$ and  $\mathcal{D}_{\mathcal{B}}$ can be defined for $D^t_k$, and Theorem~\ref{roquette plus} holds (see~\cite{dadegroup}).
\end{rem}
\section{The Roquette category} \label{section Roquette category}
\begin{mth}{Notation} \label{def X}Let $\pi$ be a projective plane over $\F_p$, and let $X$ denote a Sylow $p$-subgroup of $\Aut(\pi)\cong{\rm PGL}(3,\F_p)$. Let $\mathbb{L}$ be the set of lines of $\pi$, and let $\mathbb{P}$ be the set of points of $\pi$, both viewed as elements of $B(X)$. Let $\delta=\mathbb{L}-\mathbb{P}\in B(X)$. Equivalently
$$\delta=(X/I-X/IZ)-(X/J-X/JZ)\mvirg$$
where $I$ and $J$ are non-conjugate non-central subgroups of order $p$ of $X$, and $Z$ is the center of $X$.\par
Let $B_\delta$ denote the $p$-biset subfunctor of $B$ generated by $\delta$.
\end{mth}
\begin{rem}{Remark} When $p=2$, the group $X$ is dihedral of order 8, and $\delta$ is well-defined up to a sign. When $p>2$, the group $X$ is an extraspecial $p$-group of order $p^3$ and exponent $p$, and there are several possible choices for the element $\delta$. However, in any case, the functor $B_\delta$ does not depend on the choice of $\delta$.
\end{rem} 
\begin{mth}{Definition} The {\em Roquette category} $\mathcal{R}_p$ of finite $p$-groups is defined as the {\em idempotent additive completion} of the category $\mathcal{R}_p^\sharp$, quotient of the biset category $\mathcal{C}_p$,  defined as follows~:
\begin{enumerate}
\item the objects of $\mathcal{R}_p^\sharp$ are the finite $p$-groups.
\item if $P$ and $Q$ are finite $p$-groups, then
$$\Hom_{\mathcal{R}_p^\sharp}(P,Q)=(B/B_\delta)(Q,P)$$
is the quotient of $B(Q\times P\op)$ by $B_\delta(Q\times P\op)$.
\item the composition in $\mathcal{R}_p^\sharp$ is induced by the composition of bisets.
\item the identity morphism of the finite $p$-group $P$ in $\mathcal{R}_p^\sharp$ is the image of $\Id_P$ in $(B/B_\delta)(P,P)$.
\end{enumerate}
\end{mth}
\begin{rem}{Remark} It was shown in~\cite{rationnel} that $\mathcal{R}_p^\sharp$ is indeed a category. It was also shown there that if $p>2$, the functor $B_\delta$ is equal to the kernel $K$ of the linearization morphism $B\to R_\Q$. It follows that in this case, for any two finite $p$-groups $P$ and $Q$
$$\Hom_{\mathcal{R}_p}(P,Q)\cong R_\Q(Q\times P\op)$$
is isomorphic to the Grothendieck group of $(\Q Q,\Q P)$-bimodules, or, equivalently by the Ritter-Segal theorem, to the Grothendieck group of the subcategory of $(\Q Q,\Q P)$-permutation bimodules. In other words in this case, the category $\mathcal{R}_p$ is the full subcategory category of the category considered by Barker in~\cite{rhetoric}, consisting of finite $p$-groups. The construction of the category $\mathcal{R}_p$ is also very similar to the construction of the category $\Q G$-Morita by Hambleton, Taylor, and Williams in~\cite{htw} (Definition 1.A.4).\par
In the case $p=2$, the situation is more complicated~: the functor $B_\delta$ is a proper subfunctor of the kernel $K$, and there is a short exact sequence
$$0\to K/B_\delta \to B/B_\delta \to R_\Q\to 0$$
of $p$-biset functors. Moreover, for each $p$-group $P$, the group $(K/B_\delta)(P)$ is a finite elementary abelian 2-group of rank equal to the number of groups $S$ in a genetic basis of $P$ for which $N_P(S)/S$ is dihedral.
\end{rem}
\begin{mth}{Lemma} The direct product $(P,Q)\mapsto P\times Q$ of $p$-groups induces a well defined symmetric monoidal structure on $\mathcal{R}_p^\sharp$.
\end{mth}
\pf Let $P$, $P'$, $Q$ and $Q'$ be finite $p$-groups. If $U$ is a finite $(P',P)$-biset and $V$ is a finite $(Q',Q)$-biset, then $U\times V$ is a $(P'\times Q',P\times Q)$-biset. This induces a bilinear map
$$\pi : B(P',P)\times B(Q',Q)\to B(P'\times Q',P\times Q)\mvirg$$
and this clearly induces a symmetric monoidal structure on the biset category~$\mathcal{C}_p$. This induces a monoidal structure on the quotient category if 
$$\pi\big(B_\delta(P',P),B(Q',Q)\big)\subseteq B_\delta(P'\times Q',P\times Q)\mpoint$$
But this is a consequence of the following~: let $X$ be as defined in Notation~\ref{def X}, let $U$ be a finite $(P',P\times X)$-set, let $D$ be an $X$-set, and $V$ be a finite $(Q',Q)$-biset. Clearly, there is an isomorphism of $(P'\times Q',P\times Q)$-sets
$$(U\times_XD)\times V\cong (U\times V)\times_XD\mvirg$$
where the right action of $X$ on $U\times V$ is defined in the obvious way
$$\forall (u,v)\in U\times V,\;\forall x\in X,\;\;(u,v)x=(ux,v)\mpoint$$
The lemma follows.\findemo
\npar Recall that the objects of the idempotent additive completion $\mathcal{R}_p$ are by definition formal finite sums $\mathop{\oplus}_{(P,e)\in\mathcal{P}}\limits(P,e)$ of pairs of the form $(P,e)$, where $P$ is a finite $p$-group, and $e$ is an idempotent in $\Hom_{\mathcal{R}_p^\sharp}(P,P)=(B/B_\delta)(P,P)$. A morphism $\varphi:\mathop{\oplus}_{(P,e)\in\mathcal{P}}\limits(P,e)\to \mathop{\oplus}_{(Q,f)\in\mathcal{Q}}\limits(Q,f)$ in $\mathcal{R}_p$ is a matrix indexed by $\mathcal{P}\times\mathcal{Q}$, where the coefficient $\varphi_{(P,e),(Q,f)}$ belongs to $f\Hom_{\mathcal{R}_p^\sharp}(P,Q)e$. The composition of morphisms is given by matrix multiplication. In particular~:
\pagebreak[3]
\begin{mth}{Definition and Notation} Let $P$ be a finite $p$-group. 
\begin{enumerate}
\item The object $(P,\Id_P)$ of $\mathcal{R}_p$ is denoted by $P$. Similarly, when $Q$ is a finite $p$-group, and $f\in B(Q,P)$, the corresponding morphism from $(P,\Id_P)$ to $(Q,\Id_P)$ in the category $\mathcal{R}_p$ is simply denoted by $f$.
\item The {\em edge} $\partial P$ of $P$ is the pair $(P,f_\un^P)$ of $\mathcal{R}_p$.
\end{enumerate}
\end{mth}
The category of additive functors from $\mathcal{R}_p$ to abelian groups is equivalent to the category of additive functors from $\mathcal{R}_p^\sharp$ to abelian groups. It was shown in~\cite{rationnel} that the latter is exactly the category of rational $p$-biset functors. If $F^\sharp$ is such a functor, then $F^\sharp$ extends to a functor $F$ on $\mathcal{R}_p$ defined as follows~: 
$$F\big(\mathop{\oplus}_{(P,e)\in\mathcal{P}}\limits(P,e)\big)=\mathop{\oplus}_{(P,e)\in\mathcal{P}}\limits F^\sharp(e)\big(F^\sharp(P)\big)\mvirg$$
with the obvious definition of $F(\varphi)$ for a morphism $\varphi$ in the category $\mathcal{R}_p$. In particular, with the above notation
$$F(\partial P)=\partial F^\sharp(P)\mpoint$$
In the sequel, we will use the same symbol for $F$ and $F^\sharp$, writing in particular $F(\partial P)=\partial F(P)$.
\pagebreak[3]
\begin{mth}{Proposition} \label{sum of edges of quotients}Let $P$ be a finite $p$-group. Then
$$P\cong\dirsum{N\normal P}\partial (P/N)$$
in the category $\mathcal{R}_p$.
\end{mth}
\pf Let 
$$a:P\to\dirsum{N\normal P}\partial (P/N)$$
be the direct sum of the morphisms induced by the elements $f_\un^{P/N}\Def_{P/N}^P$ of $B(P/N,P)$, and let 
$$b: \dirsum{N\normal P}\partial (P/N)\to P$$ 
be defined similarly from the elements $\Inf_{P/N}^Pf_\un^{P/N}$ of $B(P,P/N)$.\par
By Lemma~\ref{funG}
$$\sum_{N\normal P}\Inf_{P/N}^Pf_\un^{P/N}\Def_{P/N}^P=\sum_{N\normal P}f_N^P=\Id_P$$
in $B(P,P)$, thus $a\circ b$ is equal to the identity morphism of $P$ in $\mathcal{R}_P$. Conversely, for normal subgroups $N$ and $M$ of $P$
$$f_\un^{P/N}\Def_{P/N}^P\Inf_{P/M}^Pf_\un^{P/M}=f_\un^{P/N}\Inf_{P/NM}^{P/N}\Def_{P/NM}^{P/M}f_\un^{P/M}\mvirg$$
by Equation~\ref{definf}. This is equal to 0 if $N\neq M$, by Lemma~\ref{funG}. And if $N=M$, this is equal to $f_\un^{P/N}$. It follows that $b\circ a$ is equal to the identity morphism of $\dirsum{N\normal P}\partial (P/N)$, and this completes the proof.\findemo
\begin{mth}{Corollary} \label{edge cyclic center}If $P$ is non-trivial, with cyclic center, then
$$P\cong \partial P\oplus (P/Z)$$
in $\mathcal{R}_p$, where $Z$ is the unique central subgroup of order $p$ in $P$.
\end{mth}
\pf Indeed, if $N$ is a non-trivial normal subgroup of $P$, then $N\geq Z$. Thus
$$P\cong\partial P\oplus\dirsum{N\geq Z}\partial(P/N)\cong \partial P\oplus (P/Z)$$
in the category $\mathcal{R}_p$.\findemo
\begin{rem}{Remark} \label{factor group}More generally, let $P$ be a finite $p$-group, and let $N$ be a normal subgroup of~$P$. Since in $\mathcal{R}_p$
$$P/N\cong \dirsum{N\leq M\normal P}\partial \big((P/N)\big/(M/N)\big)\cong\dirsum{N\leq M\normal P}\partial (P/M)\mvirg $$
it follows that $P/N$ is isomorphic to a direct summand of $P$ in the category~$\mathcal{R}_p$.
\end{rem}
\begin{mth}{Theorem} \label{Roquette category}\begin{enumerate}
\item The Roquette category $\mathcal{R}_p$ is an additive tensor category.
\item Let $P$ be a finite $p$-group, and $\mathcal{B}$ be a genetic basis of $P$. Then, in the category $\mathcal{R}_p$,
$$P\cong \mathop{\oplus}_{S\in\mathcal{B}}\partial \sur{N}_P(S)\mvirg$$
where $\sur{N}_P(S)=N_P(S)/S$.
\item Let $P$ be a finite $p$-group, and $\mathcal{B}$ be a genetic basis of $P$. Then, in the category $\mathcal{R}_p$,
$$\partial P\cong \oplusb{S\in\mathcal{B}}{S\cap Z(P)=\un}\partial \sur{N}_P(S)\mpoint$$
\end{enumerate}
\end{mth}
\pf Assertion 1 results from standard results~: in particular, the tensor product of $\dirsum{(P,e)\in\mathcal{P}}(P,e)$ and $\dirsum{(Q,f)\in\mathcal{Q}}(Q,f)$ is defined by
$$\big(\dirsum{(P,e)\in\mathcal{P}}(P,e)\big)\times\big(\dirsum{(Q,f)\in\mathcal{Q}}(Q,f)\big)=\oplusb{(P,e)\in\mathcal{P}}{(Q,f)\in\mathcal{Q}}(P\times Q,e\times f)\mpoint$$
For Assertion 2, by Proposition 10.7.2 of \cite{bisetfunctors}, if $F$ is a rational $p$-biset functor, the functor $F_P$ obtained from $F$ by the Yoneda-Dress construction at $P$ is also a rational $p$-biset functor. This applies in particular to the functor $Y=B/B_\delta$, so the functor $Y_P$ is rational. Hence, if $Q$ is any finite $p$-group and $\mathcal{B}_Q$ is a genetic basis of $Q$, there are mutual inverse isomorphisms
$$\xymatrix{
Y_P(Q)\ar[r]<.5ex>^-{\mathcal{D}_Q}&*!U(0.5){\mathop{\oplus}_{S\in\mathcal{B}_Q}\limits\partial Y_P\big(\sur{N}_Q(S)\big)}\ar[l]<.5ex>^-{\mathcal{I}_Q}
}
$$
where $\mathcal{I}_Q=\mathop{\oplus}_{S\in\mathcal{B}_Q}\limits\Indinf_{\sur{N}_P(S)}^P$ and $\mathcal{D}_Q=\mathop{\oplus}_{S\in\mathcal{B}_Q}\limits f_\un^{\sur{N}_P(S)}\circ\Defres_{\sur{N}_P(S)}^P$. Thus for any $f\in Y_P(Q)$
$$f=\big(\sum_{S\in\mathcal{B}_Q}\Indinf_{\sur{N}_P(S)}^Pf_\un^{\sur{N}_P(S)}\Defres_{\sur{N}_P(S)}^P\big)\circ f\mpoint$$
Applying this to $Q=P$, $\mathcal{B}_Q=\mathcal{B}$, and $f=\Id_P$ gives that $\mathcal{I}_P\circ\mathcal{D}_P=\Id_P$. On the other hand, by Proposition~6.4.4 and Theorem~9.6.1 of~\cite{bisetfunctors}, for $S,T\in\mathcal{B}$, the composition
$$f_\un^{\sur{N}_P(S)}\Defres_{\sur{N}_P(S)}^P\circ \Indinf_{\sur{N}_P(S)}^Pf_\un^{\sur{N}_P(S)}$$
is equal to $f_\un^{\sur{N}_P(S)}$ in $B\big(\sur{N}_P(S),\sur{N}_P(S)\big)$ if $T=S$, and to 0 if $T\neq S$. It follows that $\mathcal{D}_P\circ\mathcal{I}_P$ is also equal to the identity map of the direct sum $\mathop{\oplus}_{S\in\mathcal{B}}\limits\partial Y_P\big(\sur{N}_P(S)\big)$ in the category $\mathcal{R}_p$.\par
For Assertion~3, observe that by Lemma~\ref{inflation}
$$f_\un^{\sur{N}_P(S)}\Defres_{\sur{N}_P(S)}^Pf_1^P=0$$
if $S\cap Z(P)\neq \un$. Taking opposite bisets, this gives also
$$f_1^P\Indinf_{\sur{N}_P(S)}^Pf_\un^{\sur{N}_P(S)}=0\mvirg$$
so the isomorphism of Assertion~2 restricts to an isomorphism
$$\partial P\cong \oplusb{S\in\mathcal{B}}{S\cap Z(P)=\un}\partial \sur{N}_P(S)\mvirg$$
as the diagram
$$\xymatrix{
P\ar[r]^-{\cong}&\mathop{\oplus}_{S\in\mathcal{B}}\partial \sur{N}_P(S)\\
\partial P\ar[u]^{f_\un^P}\ar[r]&*!U(0.5){\oplusb{S\in\mathcal{B}}{S\cap Z(P)=\un}\limits\partial \sur{N}_P(S)}\ar@{^{(}->}[u]
}
$$
is commutative.  \findemo
\begin{mth}{Corollary} \label{zero edge} Let $P$ be a finite $p$-group. If $Z(P)$ is non-cyclic, then $\partial P=0$ in $\mathcal{R}_p$.
\end{mth}
\pf This follows from Assertion~3~: suppose indeed that there exists a genetic subgroup $S$ of $P$ such that $S\cap Z(P)\neq \un$. Then the group $Z(P)$ maps injectively in the center of the Roquette group $\sur{N}_P(S)$, which is cyclic. Hence $Z(P)$ is cyclic.\findemo
\pagebreak[3]
\begin{rem}{Examples} \label{D8Q8}
\begin{enumerate}
\item Let $P=D_8$ be a dihedral group of order 8. Let $A$, $B$, and $C$ be the subgroups of index $2$ in $P$, and let $I$ be a non-central subgroup of order~2 in $P$. Then the set $\{P,A,B,C,I\}$ is a genetic basis of $P$, and there is an isomorphism
$$P\cong \un\oplus 4\cdot \partial C_2$$
in the category $\mathcal{R}_2$, where $4\cdot \partial C_2$ denotes the direct sum of 4 copies of $\partial C_2$~: indeed, for $S\in\{A,B,C,I\}$, the group $\sur{N}_P(S)$ is isomorphic to~$C_2$.
\item Let $P=Q_8$ be a quaternion group of order 8. Let $A$, $B$, and $C$ be the subgroups of index 2 in $P$. Then the set $\{P,A,B,C,\un\}$ is a genetic basis of $P$ (such a basis is unique in this case), and there is an isomorphism
$$P\cong \un\oplus 3\cdot \partial C_2\oplus\partial Q_8$$
in the category $\mathcal{R}_2$.
\item Let $P=(C_p)^n$ be an elementary abelian $p$-group of rank $n$. Then $P$ has a unique genetic basis, consisting of $P$, and all its subgroups of index $p$. Hence
$$P\cong \un\oplus\frac{p^n-1}{p-1}\cdot\partial C_p$$
in the category $\mathcal{R}_p$.
\end{enumerate}
\end{rem}
\section{The tensor structure} \label{the tensor structure}
\begin{mth}{Notation} Let $G$ and $H$ be groups. When $L$ is a subgroup of $G\times H$, set
\begin{eqnarray*}
p_1(L)&=&\{g\in G\mid \exists h\in H,\;(g,h)\in L\}\\
p_2(L)&=&\{h\in H\mid \exists g\in G,\;(g,h)\in L\}\\
k_1(L)&=&\{g\in G\mid (g,1)\in L\}\\
k_2(L)&=&\{h\in H\mid (1,h)\in L\}\\
\end{eqnarray*}
Recall (\cite{bisetfunctors} 2.3.18 and 2.3.21) that $k_i(L)\normal p_i(L)$, for $i\in\{1,2\}$, and that $\big(k_1(L)\times k_2(L)\big)\normal L$.  Set $q(L)=L/\big(k_1(L)\times k_2(L)\big)$, and recall that there are canonical group isomorphisms
$$q(L)\cong p_1(L)/k_1(L)\cong p_2(L)/k_2(L)\mpoint$$
\end{mth}
\pagebreak[3]
\begin{mth}{Definition} Let $G$ and $H$ be groups. A subgroup $L$ of $(G\times H)$ will be called {\em diagonal} if
$$L\cap(G\times\un)=L\cap (\un\times H)=\un\mvirg$$
i.e. equivalently, if $k_1(L)=\un$ and $k_2(L)=\un$.\par
The subgroup $L$ will be called {\em centrally diagonal} if
$$L\cap\big(Z(G)\times\un\big)=L\cap\big(\un\times Z(H)\big)=\un\mvirg$$
i.e. equivalently, if $k_1(L)\cap Z(G)=\un \;\;\hbox{and}\;\;k_2(L)\cap Z(H)=\un$.
\end{mth}
\begin{mth}{Notation} Let $G$ and $H$ be groups. When $K$ is a subgroup of $G$, and $\varphi:K\to H$ is a group homomorphism, set
$$\rDelta_\varphi(K)=\{\big(x,\varphi(x)\big)\mid x\in K\}\leq G\times H\mvirg$$
and
$$\lDelta_\varphi(K)=\{\big(\varphi(x),x\big)\mid x\in K\}\leq H\times G\mpoint$$
\end{mth}
\begin{rem}{Remark} The subgroup $L$ of $(G\times H)$ is diagonal if and only if there exists a subgroup $K\leq G$ and an {\em injective} group homomorphism $\varphi:K\hookrightarrow H$ such that $L=\rDelta_\varphi(K)$.
\end{rem}
\begin{mth}{Lemma} Let $P$ and $Q$ be $p$-groups, and let $L$ and $L'$ be genetic subgroups of $(P\times Q)$ such that $L\bizlie{P\times Q}L'$. Then $L$ is centrally diagonal in $(P\times Q)$ if and only if $L'$ is centrally diagonal in $(P\times Q)$.
\end{mth}
\pf Let $(x,y)\in (P\times Q)$ such that
$$L'^{(x,y)}\cap Z_{P\times Q}(L)\leq L\mpoint$$
Since $Z(P)\times Z(Q)\leq Z_{P\times Q}(L)$, it follows that
\begin{eqnarray*}
L'\cap \big(Z(P)\times 1\big)&=&\Big(L'\cap \big(Z(P)\times 1\big)\Big)^{(x,y)}\\
&=&L'^{(x,y)}\cap \big(Z(P)\times 1\big)\\
&=&L'^{(x,y)}\cap Z_{P\times Q}(L)\cap\big(Z(P)\times 1\big)\\
&\leq& L\cap \big(Z(P)\times 1\big)=\un\mpoint
\end{eqnarray*}
A similar argument shows that $L'\cap\big(\un\times Z(Q)\big)=\un$.\findemo
Recall from Definition~\ref{axial} that an {\em axial} subgroup of a finite group $G$ is a subgroup of a cyclic subgroup of maximal order of $G$~: 
\begin{mth}{Theorem} \label{genetic product Roquette}Let $P$ and $Q$ be non-trivial Roquette $p$-groups, let $e_P$ (resp. $e_Q$) denote the exponent of $P$ (resp. of $Q$), and let $Z_P$ (resp. $Z_Q$) denote the central subgroup of order $p$ in $P$ (resp. in $Q$). \par
\begin{enumerate}
\item Let $L$ be a centrally diagonal genetic subgroup of $(P\times Q)$. 
Then $L=\rDelta_\varphi(H)$, where $H\leq P$ and $\varphi:H\hookrightarrow Q$ is an injective group homomorphism. Moreover, either $P\cong Q\cong Q_8$, and $H=P$, or $H$ is an axial subgroup of $P$ of order $\min(e_P,e_Q)$, such that $\varphi(H)$ is an axial subgroup of $Q$. 
\item Conversely~:\begin{enumerate}
\item if $P\cong Q\cong Q_8$, let $L=L_\varphi=\rDelta_\varphi(P)$, where $\varphi:P\to Q$ is a group isomorphism. Then $L$ is a centrally diagonal genetic subgroup of $(P\times Q)$, and $N_{P\times Q}(L)/L\cong C_2$. 
\item in all other cases, let $L=L_\varphi=\rDelta_\varphi(H)$, where $H$ is an axial subgroup of $P$ of order $\min(e_P,e_Q)$, and $\varphi:H\hookrightarrow Q$ is an injective group homomorphism such that $\varphi(H)$ is an axial subgroup of $Q$. Then $L$ is a centrally diagonal genetic subgroup of $(P\times Q)$. Moreover, the isomorphism class of the group $N_{P\times Q}(L)/L$ depends only on $P$ and $Q$.
\end{enumerate}
\end{enumerate}
\end{mth}
\pf $\bullet$ Observe first that the group $(Z_P\times Z_Q)L/L$ is a central subgroup of the Roquette group $N_{P\times Q}(L)/L$, hence it is cyclic. Hence $(Z_P\times Z_Q)\cap L\neq 1$. Since both $(Z_P\times \un)\cap L$ and $(\un\times Z_Q)\cap L$ are trivial, is follows that $(Z_P\times Z_Q)\cap L$ is equal to 
\begin{equation}\label{delta}
\rDelta_\psi(Z_P)=\big\{\big(z,\psi(z)\big)\mid z\in Z_P\big\}\mvirg
\end{equation}
where $\psi:Z_P\isom  Z_Q$ is some group isomorphism. In particular $p_1(L)$ contains $Z_P$, and $p_2(L)$ contains $Z_Q$.\spn
$\bullet$ Let us prove now that $L$ is diagonal, i.e. that there exists a subgroup $H$ of~$P$ and an injective group homomorphism $\varphi:H\hookrightarrow Q$ such that
$$L=\rDelta_\varphi(H)=\big\{\big(h,\varphi(h)\big)\mid h\in H\big\}\mpoint$$ 
Otherwise, at least one of the groups $k_1(L)$ or $k_2(L)$ is non trivial. But the assumption $L\cap \big(Z(P)\times\un\big)=\un$ is equivalent to $L\cap(Z_P\times\un)=\un$, i.e. $k_1(L)\cap Z_P=\un$, and similarly, the assumption $L\cap (\un\times Z(Q))=\un$ is equivalent to $k_2(L)\cap Z_Q=\un$. But if there exists a non trivial subgroup $X$ of $P$ such that $X\cap Z_P=\un$, then $p=2$, $X$ has order 2, and $P$ is dihedral or semidihedral (Lemma~\ref{en vrac}). So if $L$ is not diagonal, then $p=2$, and at least one of $P$ or $Q$ is dihedral or semidihedral.\par
Up to exchanging $P$ and $Q$, one can assume that the group $C=k_1(L)$ is non-trivial, hence non-central of order 2 in $P$. Set $A=p_1(L)$. Since $A\leq N_P(C)=CZ_P$ (Lemma~\ref{en vrac}), it follows that $q(L)=A/C$ has order 1 or~2. \par
If $q(L)=\un$, then $A=C$, and $L=C\times D$, where $D=k_2(L)=p_2(L)$. In this case 
$$N_{P\times Q}(L)/L=\big(N_P(C)/C\big)\times \big(N_Q(D)/D\big)\cong C_2\times \big(N_Q(D)/D\big) $$
cannot be a Roquette group, since $N_Q(D)/D$ is non-trivial (as $D\cap Z_Q=\un$).\par
And if $|q(L)|=2$, then $A=CZ_P$. If $(a,b)\in N_{P\times Q}(L)$, then in particular $a\in N_P(A,C)=A$. Thus $N_{P\times Q}(L)\leq A\times Q$.
Now the group $N_P(A)$ is a proper subgroup of $P$, since $A$ is elementary abelian of rank 2, and $P$ is a Roquette group. Choose $x\in P-N_P(A)$, whence $A^x\cap A=Z_P$. If $(a,b)\in L^{(x,1)}\cap (A\times Q)$, then $a\in A^x\cap A=Z_P$, hence $(a,b)={^{(x,1)}(a,b)}\in L$. Thus $L^{(x,1)}\cap(A\times Q)\leq L$, and
\begin{eqnarray*}L^{(x,1)}\cap Z_{P\times Q}(L)&\leq& L^{(x,1)}\cap N_{P\times Q}(L)\\
&\leq& L^{(x,1)}\cap(A\times Q)\\
&\leq& L
\end{eqnarray*}
but $L^{(x,1)}\neq L$, since $A^x\neq A$. It follows that $L$ is not expansive in $(P\times Q)$, hence $L$ is not a genetic subgroup of $(P\times Q)$.\spn
$\bullet$ Hence $L$ is diagonal in $(P\times Q)$, i.e.
$$L=\rDelta_\varphi(H)=\big\{\big(h,\varphi(h)\big)\mid h\in H\big\}\mvirg$$
for some subgroup $H\geq Z_P$ of $P$ and some $\varphi:H\hookrightarrow Q$, such that $\varphi(H)\geq Z_Q$. Then
$$N_{P\times Q}(L)=\Big\{(x,y)\in N_P(H)\times N_Q\big(\varphi(H)\big)\mid \forall h\in H,\;\;\varphi({^xh})={^y\varphi(h)}\Big\}\mpoint$$
The unique central subgroup of order $p$ of the Roquette group $N_{P\times Q}(L)/L$ is equal to $Z/L$, where
\begin{equation}\label{omega1Z}
Z=(Z_P\times\un)L=(\un\times Z_Q)L\mpoint
\end{equation}
For any $(x,y)\in (P\times Q)$, saying that $L^{(x,y)}\cap Z_{P\times Q}(L)$ is contained in $L$ is equivalent to saying that the group $I=L^{(x,y)}\cap Z$ is contained in $L$. In particular, for $y=1$
\begin{eqnarray*}
I&=&L^{(x,1)}\cap (Z_P\times\un)L\\
&=&\Big\{\big(h^x,\varphi(h)\big)\mid h\in H,\;\exists z\in Z_P,\;\exists h'\in H,\; \big(h^x,\varphi(h)\big)=\big(zh',\varphi(h')\big)\Big\}\\
&=&\Big\{\big(h^x,\varphi(h)\big)\mid h\in H,\;h^{-1}h^x\in Z_P\Big\}\mvirg
\end{eqnarray*}
since $\varphi(h)=\varphi(h')$ implies $h=h'$, and since $Z_P$ is central in $P$. Denoting by $[h,x]=h^{-1}h^x$ the commutator of $h$ and $x$, it follows that $I\leq L$ if and only if
$$\forall h\in H,\;\;[h,x]\in Z_P\Longrightarrow h^x\in H,\;\big(h^x,\varphi(h)\big)=\big(h^x,\varphi(h^x)\big)\mpoint$$
In other words $[h,x]\in Z_P$ implies $h^x=h$. Thus $I\leq L$ if and only if
$$\forall h\in H,\;\;[h,x]\in Z_P\Longrightarrow [h,x]=1\mpoint$$
Equivalently $[H,x]\cap Z_P=\{1\}$, where $[H,x]$ denotes {\em the set} of commutators $[h,x]$, for $h\in H$.\par
Since $L=\rDelta_\varphi(H)$ is expansive in $(P\times Q)$, it follows that 
$$[H,x]\cap Z_P=\{1\}\Longrightarrow L^{(x,1)}=L\mpoint$$
Now $(x,1)$ normalizes $L$ if and only if $x\in C_P(H)$, i.e. if $[H,x]=\{1\}$. Hence
\begin{equation}\label{expansive diagonal}
[H,x]\cap Z_P=\{1\}\Longrightarrow [H,x]=\{1\}\mpoint
\end{equation}
$\bullet$ Let us show now that unless $H=P\cong Q\cong Q_8$, the group $H$ is an axial subgroup of $P$, and the subgroup $\varphi(H)$ is an axial subgroup of $Q$.\par
Let $X$ be a cyclic subgroup of $P$ of order $e_P$, let $x$ be a generator of $X$, and suppose that $H\nleq X$. Then in particular $P$ is not cyclic, so $p=2$, the group $X$ is a (normal) subgroup of index 2 of $P$ (Lemma~\ref{en vrac}), and $X$ is equal to its centralizer in~$P$. Moreover $|H:H\cap X|=2$ since $H\cdot X=P$. The set $[H,x]$ is equal to $\{1,x^2\}$ if $P$ is cyclic or generalized quaternion, or to $\{1,x^{2+2^{n-2}}\}$ if $P$ is semidihedral~: indeed, the image of $H$ in the group of automorphisms of $X$ has order 2, as $H\cap X$ centralizes $X$, and $H$ does not. Since $Z_P$ is generated by $x^{2^{n-2}}$, it follows that $[H,x]\cap Z_P=\{1\}$ if $n\geq 4$, i.e. if $|P|\geq 16$. But $[H,x]\neq\{1\}$, hence $L$ is not expansive in $(P\times Q)$, if $P\geq 16$.\par
So if $H\nleq X$, then $p=2$, and $P$ is non-cyclic, of order at most $8$. Hence $P\cong Q_8$. If $H\neq P$, then $H$ is cyclic, and $H\nleq X$. Thus $|H|=4=e_P$, and in particular $H$ is an axis of $P$. Since $H$ embeds into~$Q$, it follows that $|H|=\min(e_P,e_Q)$. The same argument applied to $\varphi(H)$ shows that $\varphi(H)$ is an axial subgroup of $Q$, as claimed.\par
In this case moreover, the group $Q$ cannot be isomorphic to $Q_8$~: indeed otherwise, one can assume that $P=Q$ and $L=\Delta(H)$ is the diagonal embedding. Then
$$N_{P\times P}(L)=\{(a,b)\mid a^{-1}b\in H\}\mpoint$$
The group $N_{P\times P}(L)/L$ has order 8, generated by the cyclic subgroup
$$C=\{(a,1)L\mid a\in H\}$$
of index 2, and the involution $(b,b)L$, where $b\in P-H$. Hence $N_{P\times P}(L)/L\cong D_8$ is not a Roquette group, and $L$ is not a genetic subgroup of $(P\times P)$.\par
If $H$ is non-cyclic, then $H=P$, and the same argument applied to $\varphi(H)$ shows that $Q\cong Q_8$. And indeed $\rDelta_\varphi(P)$ is a genetic subgroup of $P\times Q$~: this follows from Example~\ref{diagonal}, since the map $(x,y)\mapsto \big(x,\varphi^{-1}(y)\big)$ is a group isomorphism from $(P\times Q)$ to $(P\times P)$, sending $\rDelta_\varphi(P)$ to $\Delta(P)$. Moreover $[P,P]\leq Z(P)$, and $Z(P)$ has order~2. In particular 
\begin{equation}\label{diamond 1}
N_{P\times Q}(L)/L\cong C_2
\end{equation}
does not depend on $\varphi$, up to isomorphism. This proves part a) of Assertion~2.\spn
$\bullet$ In the remaining cases $L=\rDelta_\varphi(H)$, where $H$ is a non-trivial axial subgroup of $P$, and $\varphi:H\hookrightarrow Q$ is such that $\varphi(H)$ is an axial subgroup of~$Q$. In particular $H$ is cyclic, and non-trivial. As $H\cong\varphi(H)\leq Q$, it follows that $|H|\leq\min(e_P,e_Q)$. \par
Let $C_{e_P}$ be an axis of $P$ containing $H$, and $C_{e_Q}$ be an axis of $Q$ containing $\varphi(Q)$. Then $(C_{e_P}\times C_{e_Q})/\Delta_\varphi(H)$ is an abelian normal subgroup of the Roquette group $N_{P\times Q}(L)/L$, hence it is cyclic. Thus $L=\Delta_\varphi(H)$ is not contained in the Frattini subgroup $(C_{e_P/p}\times C_{e_Q/p})$ of  $(C_{e_P}\times C_{e_Q})$. In particular $p_1(L)=C_{e_Q}$, or $p_2(L)=C_{e_Q}$. In other words $H=C_{e_P}$, or $\varphi(H)=C_{e_Q}$, hence $|H|=\min(e_P,e_Q)$. This completes the proof of Assertion~1.\spn
$\bullet$ Now assume that at least one of the groups $P$ or $Q$ is not isomorphic to~$Q_8$. Assume also that $e_P\leq e_Q$, and let $H$ be an axis of~$P$~: then $H$ is unique if $P\not\cong Q_8$, and there are three possibilities for $H$ if $P\cong Q_8$ (Lemma~\ref{en vrac}). In any case $H\normal P$. Let $K$ denote an axial subgroup of $Q$ of order~$e_P$. Such a group is unique, except if $p=2$, $e_P=4$, and $Q\cong Q_8$ (thus $P\cong C_4$ as $P$ has exponent 4, and is not isomorphic to~$Q_8$). In any case $K\normal Q$.\par
Let $\varphi: H\isom  K$ be any group isomorphism, and set $L_\varphi=\rDelta_\varphi(H)\leq (P\times Q)$. Then $L_\varphi$ is obviously centrally  diagonal, and
$$N_{P\times Q}(L_\varphi)=\big\{(a,b)\in (P\times Q)\mid \forall h\in H,\;\;\varphi({^ah})={^b\varphi(h)}\big\}\mpoint$$
Since $H$ is cyclic of order $e_P$, the map
$$\pi_H:r\in(\Z/e_P\Z)^\times\mapsto (x\mapsto x^r)\in\Aut(H)$$
is a canonical group isomorphism. Similarly, the map
$$\pi_K:r\in(\Z/e_P\Z)^\times\mapsto (x\mapsto x^r)\in\Aut(K)$$
is a canonical group isomorphism.\par
Let $\alpha:P\to \Aut(H)\stackrel{\pi_H^{-1}}{\longrightarrow}(\Z/e_P\Z)^\times$ denote the group homomorphism obtained from the action of $P$ on its normal subgroup~$H$ by conjugation, and $\beta:Q\to \Aut(K)\stackrel{\pi_K^{-1}}{\longrightarrow}(\Z/e_P\Z)^\times$ denote the group homomorphism obtained from the action of $Q$ on its normal subgroup~$K$. Then
$$N_{P\times Q}(L_\varphi)=\big\{(a,b)\in (P\times Q)\mid \alpha(a)=\beta(b)\big\}\mpoint$$
Now the group $(\Z/e_P\Z)^\times$ is abelian. The map 
$$\Theta: (a,b)\in (P\times Q)\mapsto \beta(b)^{-1}\cdot\alpha(a)\in(\Z/e_P\Z)^\times$$
is a group homomorphism, and $N_{P\times Q}(L_\varphi)=\Ker\;\Theta$. In particular, it is a normal subgroup of $(P\times Q)$, which does not depend on $\varphi$, once $H$ and $K=\varphi(H)$ are fixed. In particular $L_\varphi$ is an expansive subgroup of $(P\times Q)$, by Example~\ref{normalisateur normal}.
Moreover setting $I_{P,Q}=\Im(\alpha)\cap\Im(\beta)$, there is an exact sequence
\begin{equation}\label{simple ses}
\un\to C_P(H)\times C_Q(K)\to N_{P\times Q}(L_\varphi)\stackrel{\Psi}{\longrightarrow} I_{P,Q}\to\un\mvirg
\end{equation}
where $\Psi(a,b)=\alpha(a)=\beta(b)$, for $(a,b)\in N_{P\times Q}(L_\varphi)$.\mmp
\S.~Suppose first that $e_P=p$, i.e. that $P\cong C_p$. In this case $H=Z_P=P$, and $K=Z_Q$, so $L_\varphi$ is central in $(P\times Q)$. Moreover
\begin{equation}\label{diamond 2}
N_{P\times Q}(L_\varphi)/L_\varphi=(P\times Q)/L_\varphi=(P\times Q)/\rDelta_\varphi(P)\cong Q
\end{equation}
is a Roquette group, independent of $\varphi$, up to isomorphism. In particular $L_\varphi$ is a genetic subgroup of $(P\times Q)$.\mmp
\S.~Assume from now on that $e_P\geq p^2$. Then $C_P(H)\cong C_{e_P}$, and $C_Q(K)\cong C_{e_Q}$, by Lemma~\ref{en vrac}.\mmp 
\S\S.~If $p>2$, then $P$ and $Q$ are cyclic, hence $H=P$, and 
\begin{equation}\label{diamond 3}N_{P\times Q}(L_\varphi)/L_\varphi\cong (P\times Q)/\rDelta_\varphi(P)\cong Q
\end{equation}
as above. It is a Roquette group, independent of $\varphi$, up to isomorphism. In particular $L_\varphi$ is genetic in $(P\times Q)$.
\mmp
\S\S.~Assume now that $p=2$. The image of $\alpha$ has order $|P:C_P(H)|$, which is equal to 1 if $P$ is cyclic, and to 2 otherwise. Similarly, the image of $\beta$ has order $|Q:C_Q(K)|$, which is equal to 1 if $K$ is central in $Q$, i.e. if $Q$ is cyclic (since $|K|=e_P\geq 4$ by assumption), and to 2 otherwise. Set $I_{P,Q}=\Im(\alpha)\cap\Im(\beta)$. Then $I_{P,Q}$ has order 1 or 2, and there is an exact sequence
\begin{equation}\label{type}
1\to C_{e_P}\times C_{e_Q}\to N_{P\times Q}(L_\varphi)\to I_{P,Q}\to 1\mpoint
\end{equation}
Note that $I_{P,Q}$ does not depend on $\varphi:H\isom K$~: more precisely 
$$\Im(\alpha)=\left\{\begin{array}{cl}\{1\}&\hbox{if $P$ is cyclic,}\\
\{1,-1\}&\hbox{if $P$ is dihedral or generalized quaternion,}\\
\{1,e_P/2-1\}&\hbox{if $P$ is semidihedral.}
\end{array}\right.$$
So $\Im(\alpha)$ only depends on the type of $P$.\par
Similarly
$$\Im(\beta)=\left\{\begin{array}{cl}\{1\}&\hbox{if $Q$ is cyclic,}\\
\{1,-1\}&\hbox{if $Q$ is dihedral or generalized quaternion,}\\
\{1,e_Q/2-1\}&\hbox{if $Q$ is semidihedral.}
\end{array}\right.$$
Moreover if $Q$ is semidihedral and if $e_Q>e_P$, then $e_Q/2-1\equiv -1 ({\rm mod.} e_P)$, hence $\Im(\beta)=\{1,-1\}$. In other words, the group $I_{P,Q}$ is trivial in one of the following cases~:
$$\!\!\left\Arrowvert\begin{array}{l}
\bullet\;\hbox{ $P$ or $Q$ is cyclic,}\\
\bullet\;\hbox{ $P$ is dihedral or generalized quaternion, $Q$ is semidihedral, and $e_P=e_Q$,} \\
\phantom{\bullet}\;\hbox{(i.e. equivalently $|P|=|Q|$),}\\
\bullet \;\hbox{ $P$ is semidihedral, and $Q$ is dihedral or generalized quaternion,}\\
\bullet\;\hbox{ $P$ and $Q$ are semidihedral, and $e_P<e_Q$ (i.e. equivalently $|P|<|Q|$),}
\end{array}\right.$$
and the group $I_{P,Q}$ has order 2 in all other cases, i.e. in one of the following cases~: 
$$\!\!\left\Arrowvert\begin{array}{l}
\bullet\;\hbox{ $P$ and $Q$ are dihedral or generalized quaternion,}\\
\bullet\;\hbox{ $P$ is dihedral or generalized quaternion, $Q$ is semidihedral, and $|Q|{>}|P|$,}\\
\bullet\;\hbox{ $P$ and $Q$ are semidihedral, and $P\cong Q$.}
\end{array}\right.$$
As $L_\varphi\leq C_{e_P}\times C_{e_Q}$, the exact sequence~\ref{type} yields the exact sequence
\begin{equation}\label{exact diamond}
1\to (C_{e_P}\times C_{e_Q})/L_\varphi\to N_{P\times Q}(L_\varphi)/L_\varphi\to I_{P,Q}\to 1\mpoint
\end{equation}
\sou{Case 1~:} If $I_{P,Q}$ is trivial, then $N_{P\times Q}(L_\varphi)\cong C_{e_P}\times C_{e_Q}$, and 
\begin{equation}\label{diamond 4}
N_{P\times Q}(L_\varphi)/L_\varphi\cong C_{e_Q}\mvirg
\end{equation}
which is a Roquette group, independent of $\varphi$, up to isomorphism. In particular $L_\varphi$ is a genetic subgroup of $(P\times Q)$. \spn
\sou{Case 2~:} Suppose now that $I_{P,Q}$ has order 2, i.e. that $\Im(\alpha)=\Im(\beta)=\{1,\epsilon\}$, where $\epsilon$ is either $-1$ or $e_Q/2-1$ (in the case where $P$ and $Q$ are semidihedral and isomorphic). One can choose an element $u\in P$, of order 2 if $P$ is dihedral or semidihedral, and of order 4 if $P$ is generalized quaternion, such that $\alpha(u)=\epsilon$. Similarly, one can choose an element $v\in Q$, of order 2 if $Q$ is dihedral or semidihedral, and of order 4 if $Q$ is generalized quaternion, such that $\beta(v)=\epsilon$. These choices imply that $(u,v)\in N_{P\times Q}(L_\varphi)$.\par
In the exact sequence~\ref{exact diamond}
$$
1\to (C_{e_P}\times C_{e_Q})/L_\varphi\to N_{P\times Q}(L_\varphi)/L_\varphi\to \{1,\epsilon\}\to 1\mvirg
$$
the group $C=(C_{e_P}\times C_{e_Q})/L_\varphi$ is cyclic, isomorphic to $C_{e_Q}$. The element $\pi=(u,v)L_\varphi$ of $N_{P\times Q}(L_\varphi)/L_\varphi$ acts on $C$ in the same way that $v$ acts on the subgroup $C_{e_Q}$ of $Q$, namely by inversion if $Q$ is dihedral or generalized quaternion, and by raising elements to the power $e_Q/2-1$ if $Q$ is semidihedral. Finally $\pi^2=(u^2,v^2)L_\varphi=L_\varphi$ if none of $P$ and $Q$ are generalized quaternion. If $P$ is generalized quaternion and $Q$ is not, then $\pi^2=(z_P,1)L_\varphi\in C-\{1\}$, where $z_P$ is a generator of $Z_P$. Similarly, if $Q$ is generalized quaternion and $P$ is not, then $\pi^2=(1,z_Q)L_\varphi\in C-\{1\}$, where $z_Q$ is a generator of $Z_Q$. In these two cases $\pi^4=(1,1)L_\varphi=L_\varphi$, so $\pi$ has order 4 in $N_{P\times Q}(L_\varphi)/L_\varphi$. And finally, if both $P$ and $Q$ are generalized quaternion, then $\pi^2=(z_P,z_Q)L_\varphi=L_\varphi$, since $\varphi(Z_P)=Z_Q$. \par
\begin{equation}\label{diamond 5}
\hbox{It follows that the group $N_{P\times Q}(L_\varphi)/L_\varphi$ has order $2e_Q=|Q|$, and that it is~:}
\end{equation}
$$
\!\!\left\Arrowvert\begin{array}{l} 
\bullet\;\hbox{ dihedral if $P$ and $Q$ are both dihedral, or both generalized quaternion,}\\
\bullet\;\hbox{ generalized quaternion if one of $P$, $Q$ is generalized quaternion, and the}\\
\phantom{\bullet}\;\hbox{ other is dihedral,}\\
\bullet\;\hbox{ semidihedral if $Q$ is semidihedral.}
\end{array}\right.
$$
So $N_{P\times Q}(L_\varphi)/L_\varphi$ is a Roquette group, independent of $\varphi$, up to isomorphism. In particular, it follows that $L_\varphi$ is a genetic subgroup of $(P\times Q)$. This completes the proof of Theorem~\ref{genetic product Roquette}.\findemo
\begin{mth}{Notation}\label{diamond}\begin{itemize}
\item Let $P$ and $Q$ be Roquette $p$-groups. If $P$ and $Q$ are non-trivial, set $P\diamond Q=N_{P\times Q}(L)/L$, where $L$ is a centrally diagonal genetic subgroup of $(P\times Q)$. Set moreover $\un\diamond P=P\diamond\un=P$.
\item Let $P$ and $Q$ be Roquette $p$-groups. If $P$ and $Q$ are non-trivial, let $\nu_{P,Q}$ denote the number of equivalence classes of centrally diagonal genetic subgroups of $(P\times Q)$ for the relation $\bizlie{P\times Q}$. Set moreover $\nu_{\un,P}=\nu_{P,\un}=1$.
\end{itemize}
\end{mth}
\vspace{1ex}
\begin{rem}{Remark} if $P=\un$, and $Q$ is a Roquette $p$-group, then $P\times Q\cong Q$, and the centrally diagonal genetic subgroups of $P\times Q$ are the subgroups $\un\times R$, where $R$ is a genetic subgroup of $Q$ such that $R\cap Z(Q)=\un$. The only such subgroup is $R=\un$, so $N_Q(R)/R\cong Q\cong N_{P\times Q}(\un\times R)/(\un\times R)$. Hence the above definition of $P\diamond Q$ and $\nu_{P,Q}$ is consistent, in the case $P=\un$.
\end{rem}
\pagebreak[3]
\begin{mth}{Theorem}\label{structure diamond} Let $P$ and $Q$ be Roquette $p$-groups, of exponents $e_P$ and $e_Q$, respectively. Suppose $e_P\leq e_Q$, and set $q=|Q|$. Then
$$P\diamond Q\cong\left\{\begin{array}{cl}
Q&\hbox{if $P=\un$ or $P\cong C_p$,}\\
C_2&\hbox{if $P\cong Q\cong Q_8$,}\\
D_q&\hbox{if $q\geq 16$ and $P$ and $Q$ are both dihedral,}\\
&\hbox{\ \ \ or both generalized quaternion,}\\
Q_q&\hbox{if one of $P, Q$ is dihedral,}\\
&\hbox{\ \ \ and the other one is generalized quaternion,}\\
SD_q&\hbox{if $Q$ is semidihedral, and}\\
&\hbox{- either $P$ is dihedral or generalized quaternion,}\\
&\hspace{\parindent}\hbox{and $|P|<|Q|$,}\\
&\hbox{- or $P\cong Q$,}\\
C_{e_Q}&\hbox{otherwise.}
\end{array}\right.$$
\end{mth}
\pf The case $P=\un$ is trivial, the case $P\cong C_p$ follows from~\ref{diamond 2}, the case $P\cong Q\cong Q_8$ follows from~\ref{diamond 1}, the three next cases in the list follow from~\ref{diamond 5}, and the last case follows from~\ref{diamond 3} and~\ref{diamond 4}.\findemo
\begin{mth}{Theorem} \label{tensor}Let $P$ and $Q$ be Roquette $p$-groups, of exponent $e_P$ and~$e_Q$, respectively, and let $m=\min(e_P,e_Q)$. \begin{enumerate}
\item If $p=2$ and and one of the groups $P$ or $Q$ is isomorphic to $Q_8$, then then $L\bizlie{P\times Q}L'$ for any centrally diagonal genetic subgroups $L$ and $L'$ of $P\times Q$. In other words $\nu_{P,Q}=1$.
\item In all other cases, if $L$ and $L'$ are centrally diagonal genetic subgroups of $P\times Q$, then $L\bizlie{P\times Q}L'$ if and only if $L$ and $L'$ are conjugate in $P\times Q$. In particular $\nu_{P,Q}=\phi(m)m\dsp\frac{|P\diamond Q|}{|P||Q|}$, where $\phi$ is the Euler function. 
\end{enumerate}
\end{mth}
\pf Assume $e_P\leq e_Q$, without loss of generality.\mpn
1) If $P=\un$, then $\nu_{P,Q}=1$, and $P\diamond Q=Q$ by definition, and $e_P=1$, so there is nothing to prove.\mpn
2) If $p=2$ and $P\cong Q\cong Q_8$, then by Theorem~\ref{genetic product Roquette}, a genetic centrally diagonal subgroup~$L$ of $P\times Q$ is of the form $L_\varphi=\rDelta_\varphi(P)$, where $\varphi:P\to Q$ is some group isomorphism. Moreover $N_{P\times Q}(L)/L\cong C_2$, so $P\diamond Q\cong C_2$.\par
Now let $\varphi,\psi:P\to Q$ be two group isomorphisms. Then $L_\varphi\bizlie{P\times Q}L_\psi$ if and only if there exists $(x,y)\in (P\times Q)$ such that
\begin{eqnarray*}
L_\varphi^{(x,y)}\cap Z_{P\times Q}(L_\psi)&\leq& L_\psi\\
^{(x,y)}L_\psi\cap Z_{P\times Q}(L_\varphi)&\leq& L_\varphi\mpoint
\end{eqnarray*}
These conditions depend only on the double coset $N_{P\times Q}(L_\varphi)(x,y)N_{P\times Q}(L_\psi)$, which admits a representative of the form $(u,1)$. \par
Now the condition $L_\varphi^{(u,1)}\cap Z_{P\times Q}(L_\psi)\leq L_\psi$ is equivalent to
$$\forall h\in P,\;\;\varphi(h)^{-1}\psi(h^u)\in Z_Q \Longrightarrow \varphi(h)^{-1}\psi(h^u)=1\mvirg$$
and similarly, the condition $^{(u,1)}L_\psi\cap Z_{P\times Q}(L_\varphi)\leq L_\varphi$ is equivalent to
$$\forall h\in P,\;\;\psi(h)^{-1}\varphi({^uh})\in Z_Q \Longrightarrow \psi(h)^{-1}\varphi({^uh})=1\mpoint$$
Applying $\psi^{-1}$ to the first condition and $\varphi^{-1}$ to the second one, and setting $\theta=\psi^{-1}\varphi$, these two conditions become
\begin{eqnarray*}
\forall h\in P,\;\;\theta(h)^{-1}h^u\in Z_P&\Longrightarrow& \theta(h)=h^u\\
\forall h\in P,\;\;\theta^{-1}(h)^{-1}{^uh}\in Z_P&\Longrightarrow&\theta^{-1}(h)={^uh}\mpoint
\end{eqnarray*}
Since $h^uh^{-1}\in[P,P]=Z_P$, there are equivalences
\begin{eqnarray*}
\theta(h)^{-1}h^u\in Z_P\;\; \Longleftrightarrow&\theta(h)^{-1}h\in Z_P&\Longleftrightarrow\;\; h^{-1}\theta(h)\in Z_P\mvirg\\
\theta^{-1}(h)^{-1}{^uh}\in Z_P\;\;\Longleftrightarrow&\theta^{-1}(h)^{-1}{h}\in Z_P&
\Longleftrightarrow\;\; h^{-1}\theta(h)\in Z_P\mpoint
\end{eqnarray*}
Hence, in order to prove that $L_\varphi\bizlie{P\times Q}L_\psi$ for any $\varphi, \psi:P\isom Q$, it is enough to prove that 
\begin{equation}\label{condautQ8}\forall \theta\in \Aut(P),\;\exists u\in P,\;\forall h\in P,\;\;\theta(h)h^{-1}\in Z_P\implies \left\{\begin{array}{l}\theta(h)=h^u\\\theta^{-1}(h)={^uh}\end{array}\right.\mpoint
\end{equation}
If $h$ has order 1 or 2, then $\theta(h)=h=h^u$ for any $u\in P$, hence the conditions $\theta(h)=h^u$ and $\theta^{-1}(h)={^uh}$ only have to be checked for $|h|=4$. Now the group $\Aut(P)$ permutes the three cyclic subgroups of order 4 of $P$, and this gives an exact sequence
$$1\to\Inn(P)\to \Aut(P)\to S_3\to 1\mvirg$$
where $S_3$ is the symmetric group on three symbols. Saying that $\theta(h)h^{-1}\in Z_P$ is equivalent to saying that $\theta({<}h{>})={<}h{>}$. Hence, either $\theta$ stabilizes the three subgroups of order 4 of $P$, and in this case $\theta$ is inner, hence there exists $u\in P$ such that $\theta(h)=h^u$, for any $h\in P$, hence $\theta^{-1}(h)={^uh}$, for any $h\in P$. \par
Or there exists a unique subgroup $C$ of order $4$ of $P$ such that $\theta(C)=C$. Then either $\theta(h)=h$ for any $h\in C$, or $\theta(h)=h^{-1}$, for any $h\in C$. In the first case, take $u=1$, and in the second case take $u\in P-{<}h{>}$, and then $\theta(h)=h^u$ for any $h\in C$, hence $h=\theta^{-1}(h)^u$ for $h\in C$, since $C=\theta(C)$. Hence~\ref{condautQ8} holds. This completes the proof in this case.\mpn
3) If $P\cong Q_8$ and $Q\ncong Q_8$, then a centrally diagonal genetic subgroup of $(P\times Q)$ is of the form $L=\rDelta_\varphi(H)$, where $H$ is one of the 3 subgroups of order 4 of $P$, and $\varphi$ is some isomorphism from $H$ to the unique axial subgroup $K$ of order 4 of $Q$. Moreover $Z_{P\times Q}(L)=L(\un\times Z_Q)$.\par
Let $H$ and $H'$ be subgroups of order 4 of $P$. Let $\varphi: H\to K$ and $\varphi':H'\to\nolinebreak K$ be group isomorphisms, and set $L=\rDelta(\varphi)$ and $L'=\rDelta_{\varphi'}(H')$. Suppose first that $H\neq H'$, and let $(a,b)\in L'\cap Z_{P\times Q}(L)=L(\un\times Z_Q)$. It means that $a\in H'\cap H=Z_P$, and that there exists $z\in Z_Q$ such that $\varphi'(a)=\varphi(a)z$. But the restrictions of $\varphi$ and $\varphi'$ to $Z_P$ are equal, so $\varphi'(a)=\varphi(a)$, hence $z=1$. It follows that $L'\cap Z_{P\times Q}(L)\leq L$, hence $L\cap Z_{P\times Q}(L')\leq L'$ by symmetry, so $L\bizlie{P\times Q} L'$ in this case.\par
Now if $H=H'$, choose a subgroup $H''$ of order 4 in $P$, different from~$H$, and a group isomorphism $\varphi'':H''\to K$. Set $L''=\rDelta_{\varphi''}(H'')$. Then $L\bizlie{P\times Q}L''\bizlie{P\times Q}L'$, by the previous argument, thus $L\bizlie{P\times Q}L'$. \par
Hence $\nu_{P,Q}=1$ in this case, as was to be shown.\mpn
4) In the case there are several choices for $K$, i.e. if $Q\cong Q_8$ and $K$ has order $4=\min(e_P,e_Q)$, it follows that $P\cong C_4$, since $P\ncong Q_8$. In this case, we can exchange $P$ and $Q$, and use the previous argument. Hence $\nu_{P,Q}=1$ in this case as well.\mpn
5) In all other cases, by Theorem~\ref{genetic product Roquette}, a centrally diagonal genetic subgroup $L$ of $(P\times Q)$ is of the form $\rDelta_\varphi(H)$, where $H\leq P$ is the unique axis of $P$, and $\varphi:H\hookrightarrow Q$ is a group isomorphism to the unique axial subgroup $K$ of order $e_P$ of $Q$. The normalizer of $L$ in $P\times Q$ does not depend on $\varphi$, by~\ref{simple ses}, so it does not depend on $L$, since $H$ and $K$ are also unique.\par
Let $L$ and $L'$ be two such centrally diagonal genetic subgroups of $P\times Q$. Then $N_{P\times Q}(L)=N_{P\times Q}(L')$, thus $L\bizlie{P\times Q}L'$ if and only if $L$ and $L'$ are conjugate in $(P\times Q)$, by Lemma~\ref{normalisateur normal bizlie}. Moreover, it follows from the definition of $P\diamond Q$ that
$$|N_{P\times Q}(L)=|L||P\diamond Q|=e_P|P\diamond Q|\mvirg$$
so the conjucacy class of $L$ in $(P\times Q)$ has cardinality $\dsp\frac{|P||Q|}{e_P|P\diamond Q|}$. Since there are $\phi(e_P)$ possible choices for the isomorphism $\varphi:H\to K$, i.e. $\phi(e_P)$ centrally diagonal subgroups of $(P\times Q)$, it follows that
$$\nu_{P,Q}=\phi(e_P)e_P\frac{|P\diamond Q|}{|P||Q|}\mvirg$$
as was to be shown.\findemo
\begin{rem}{Remark} Suppose that $P\cong Q_8$ and $Q\ncong Q_8$. Then $|P\diamond Q|=|Q|$, by Theorem~\ref{structure diamond}. Hence 
$$\phi(e_P)e_P\frac{|P\diamond Q|}{|P||Q|}=2\times 4\times\frac{|Q|}{8|Q|}=1=\nu_{P,Q}\mvirg$$
so the formula for $\nu_{P,Q}$ holds in this case. The only case where $\nu_{P,Q}$ is not equal to $\phi(m)m\dsp\frac{|P\diamond Q|}{|P||Q|}$ (where $m=\min(e_P,e_Q)$) is the case $P\cong Q\cong Q_8$~: in this case $\nu_{P,Q}=1$, but
$$\phi(m)m\frac{|P\diamond Q|}{|P||Q|}=2\times 4\times\frac{2}{8\times 8}=\frac{1}{4}\mpoint$$
\end{rem}
\pagebreak[3]
\begin{mth}{Corollary} \label{edge product}Let $P$ and $Q$ be Roquette $p$-groups. Then, in the category $\mathcal{R}_p$
$$\partial P\times\partial Q\cong \nu_{P,Q}\cdot\partial(P\diamond Q)\mpoint$$
In other words, if $P$ has exponent $e_P$, if $Q$ has order $q$ and exponent $e_Q$, and if $e_P\leq e_Q$~:
$$\partial P\times\partial Q=\left\{\begin{array}{ll}\partial Q&\hbox{if $P=\un$ or $P\cong C_2$,}\\
\rule{0ex}{4ex}\partial C_2&\hbox{if $P\cong Q\cong Q_8$,}\\
\rule{0ex}{4ex}\dsp\frac{\phi(e_P)}{2}\cdot\partial D_q&\hbox{if $q\geq 16$ and $P$ and $Q$ are both dihedral,}\\
&\hbox{\ \ \ or both generalized quaternion,}\\
\dsp\frac{\phi(e_P)}{2}\cdot\partial Q_q&\hbox{if one of $P, Q$ is generalized quaternion,}\\
&\hbox{\ \ \  and the other one is dihedral,}\\
\dsp\frac{\phi(e_P)}{2}\cdot\partial SD_q&\hbox{if $Q$ is semidihedral, and}\\
&\hbox{- either $P$ is dihedral or generalized}\\
&\hspace{\parindent}\hbox{quaternion, and $|P|<|Q|$,}\\
&\hbox{- or $P\cong Q$,}\\
\dsp\rule{0ex}{4ex}\frac{\phi(e_P)e_Pe_Q}{|P||Q|}\cdot\partial C_{e_Q}&\hbox{otherwise.}   
\end{array}\right.$$
\end{mth}
\pf Let $\mathcal{B}$ be a genetic basis of the group $R=P\times Q$. In the category $\mathcal{R}_p$, the product $\partial P\times\partial Q$ is equal to $(P\times Q,f_\un^P\times f_\un^Q)$, and it is a summand of $R$. By Theorem~\ref{Roquette category}, there are mutual inverse isomorphisms
$$\xymatrix{
R\ar[r]<.5ex>^-{\mathcal{D}}&*!U(0.4){\mathop{\oplus}_{S\in \mathcal{B}}\limits\partial\sur{N}_{R}(S)}\ar[l]<.5ex>^-{\mathcal{I}}\mvirg
}
$$
where $\mathcal{I}$ is the direct sum of the maps $\Indinf_{\sur{N}_{R}(S)}^{R}$, and $\mathcal{D}$ is the direct sum of the maps $f_\un^{\sur{N}_{R}(S)}\Defres_{\sur{N}_{R}(S)}^{R}$. Corollary~\ref{edge product} follows from the fact that
\begin{equation}\label{restricted iso}f_\un^{\sur{N}_R(S)}\Defres_{\sur{N}_{R}(S)}^{R}(f_\un^P\times f_\un^Q)=0
\end{equation}
unless $S$ is centrally diagonal in $R=P\times Q$~: indeed $\Defres_{\sur{N}_{R}(S)}^{R}$ is given by the $\big(\sur{N}_{R}(S),R\big)$-biset $S\dom R$. On the other hand
\begin{eqnarray*}
f_\un^P\times f_\un^Q&=&(P/\un-P/Z_P)\times (Q/\un-Q/Z_Q)\\
&=&R/(\un\times \un)-R/(\un\times Z_Q)-R/(Z_P\times \un)+R/(Z_P\times Z_Q)\mvirg
\end{eqnarray*}
hence $\Defres_{\sur{N}_{R}(S)}^{R}(f_\un^P\times f_\un^Q)$ is equal to
$$S\dom R/(\un\times \un)-S\dom R/(\un\times Z_Q)-S\dom R/(Z_P\times \un)+S\dom R/(Z_P\times Z_Q)$$
which is
\begin{equation}\label{zero}
S\dom  R-S(\un\times Z_Q)\dom R-S(Z_P\times \un)\dom  R+S(Z_P\times Z_Q)\dom R\mpoint
\end{equation}
If $S$ is not centrally diagonal in $R=P\times Q$, then either $S=S(\un\times Z_Q)$ or $S=S(Z_P\times\un)$. In each case the sum~\ref{zero} vanishes.\par
And if $S$ is centrally diagonal in $R$, then 
$$S(\un\times Z_Q)=S(Z_P\times \un)=S(Z_P\times Z_Q)\mvirg$$
since the image of these groups in the Roquette group $\sur{N}_R(S)$ is equal to its unique central subgroup $\widehat{S}/S$ of order $p$. In this case
$$\Defres_{\sur{N}_{R}(S)}^{R}(f_\un^P\times f_\un^Q)=S\dom R-\widehat{S}\dom R\mpoint$$
Since $f_\un^{\sur{N}_R(S)}=N_R(S)/S-N_R(S)/\widehat{S}$, it follows that $\Defres_{\sur{N}_{R}(S)}^{R}(f_\un^P\times f_\un^Q)$ is invariant by composition with $f_\un^{\sur{N}_R(S)}$. \par
Conversely, if $S$ is not centrally diagonal in $(P\times Q)$, then 
$$(f_\un^P\times f_\un^Q)\Indinf_{\sur{N}_R(S)}^Rf_\un^{\sur{N}_R(S)}=0\mvirg$$
as can be seen by taking opposite bisets in equation~\ref{restricted iso}. And if $S$ is centrally diagonal, then
$$(f_\un^P\times f_\un^Q)\Indinf_{\sur{N}_R(S)}^Rf_\un^{\sur{N}_R(S)}=R/S-R/\widehat{S}=\Indinf_{\sur{N}_R(S)}^Rf_\un^{\sur{N}_R(S)}\mpoint$$
Hence the isomorphisms $\mathcal{D}$ and $\mathcal{I}$ restrict to mutual inverse isomorphisms between $\partial P\times \partial Q$ and the direct sum of the edges $\partial \sur{N}_R(S)$, where $S$ is a centrally diagonal genetic subgroup of $R$. But for all such subgroups $S$, the group $\sur{N}_R(S)$ is isomorphic to $P\diamond Q$, and there are $\nu_{P,Q}$ centrally diagonal subgroups in a genetic basis of $R=P\times Q$. This completes the proof.\findemo
\vspace{-2ex}
\section{Examples and applications} \label{examples and applications}
\vspace{-1ex}
\npar Suppose first that $p$ is odd. Then the Roquette $p$-groups are just the cyclic groups $C_{p^n}$, for $n\geq 0$. The ``multiplication rule" of the edges $\partial C_{p^n}$ is the following
\begin{equation}\label{product cyclic edges}
\forall m,\forall n\in \N,\;\partial C_{p^m}\times \partial C_{p^n}=\phi(p^{\min(m,n)})\partial C_{p^{\max(m,n)}}\mvirg
\end{equation}
where $\phi$ is the Euler function (thus $\phi(p^k)=p^{k-1}(p-1)$ if $k>0$, and $\phi(1)=1$).
\pagebreak[3]
\npar Some surprising phenomenons occur when $p=2$~:
\begin{mth}{Proposition}\label{edge C2}
In $\mathcal{R}_2$, the edge $\partial C_2$ is isomorphic to the trivial group~$\un$ (or its edge $\partial \un$). 
\end{mth}
\pf Indeed Corollary~\ref{zero edge} implies that if $E\cong (C_2)^2$, then $\partial E=0$ in $\mathcal{R}_2$. Let $X$, $Y$ and $Z$ denote the subgroups of order 2 of $E$. The element
$$u=\Res_X^E\times_Ef_\un^E\times_E\Ind_Y^E$$
of $B(X,Y)$ can be viewed as a morphism from $Y$ to $X$ in the category $\mathcal{R}_2$, which factors through $\partial E$. So this morphism is equal to 0. Since 
$$f_\un^E=E/\un-E/X-E/Y-E/Z+2E/E\mvirg$$
it follows that
$$u=\Ind_\un^X\Res_\un^Y-\Inf_\un^X\Res_\un^Y-\Ind_\un^X\Def_\un^Y-\Iso(\varphi)+2\Inf_\un^X\Def_\un^Y\mvirg$$
where $\varphi$ is the unique group isomorphism from $Y$ to $X$.
Thus
$$0=u\;\Iso(\varphi^{-1})=\Ind_\un^X\Res_\un^X-\Inf_\un^X\Res_\un^X-\Ind_\un^X\Def_\un^X-\Id_X+2\Inf_\un^X\Def_\un^X\mpoint$$
Hence in the category $\mathcal{R}_2$
$$\Id_X=(\Ind_\un^X-\Inf_\un^X)(\Res_\un^X-\Def_\un^X)+\Inf_\un^X\Def_\un^X\mpoint$$
It follows that
\begin{equation}\label{ab}f_\un^X=f_\un^X(\Ind_\un^X-\Inf_\un^X)(\Res_\un^X-\Def_\un^X)f_\un^X\mpoint
\end{equation}
But on the other hand $f_\un^X=\Id_X-\Inf_\un^X\Def_\un^X$, so
\begin{eqnarray*}
f_\un^X(\Ind_\un^X-\Inf_\un^X)&=&f_\un^X\Ind_\un^X\\
&=&\Ind_\un^X-\Inf_\un^X\Def_\un^X\Ind_\un^X\\
&=&\Ind_\un^X-\Inf_\un^X\mpoint
\end{eqnarray*}
It follows that
\begin{eqnarray*}
(\Res_\un^X-\Def_\un^X)f_\un^X(\Ind_\un^X-\Inf_\un^X)&=&(\Res_\un^X-\Def_\un^X)(\Ind_\un^X-\Inf_\un^X)\\
&=&2\Id_\un-\Id_\un-\Id_\un+\Id_\un\\
&=&\Id_\un\mpoint
\end{eqnarray*}
Thus, setting 
\begin{eqnarray*}
a&=&f_\un^X(\Ind_\un^X-\Inf_\un^X)\in\Hom_{\mathcal{R}_2}(\un,\partial X)\\
b&=&(\Res_\un^X-\Def_\un^X)f_\un^X\in\Hom_{\mathcal{R}_2}(\partial X,\un)\mvirg
\end{eqnarray*}
the composition $b\circ a$ is equal to $\Id_\un$, and Equation~\ref{ab} shows that the composition $a\circ b$ is equal to the identity of~$\partial X$. So $a$ and $b$ are mutual inverse isomorphisms between $\un$ and $\partial X$.\findemo
\begin{mth}{Corollary} Let $F$ be a rational $2$-biset functor. Then for any finite 2-group $P$
$$F(C_2\times P)\cong F(P)\oplus F(P)\mpoint$$
$$F(D_8\times P)\cong F(P)^{\oplus 5}\mpoint$$
\end{mth}
\pf Indeed rational $p$-biset functors are exactly those $p$-biset functors which factor through the category $\mathcal{R}_p$. And in the category $\mathcal{R}_2$, by Theorem~\ref{Roquette category}, there is an isomorphism
$$C_2\cong \un\oplus \partial C_2\cong \un\oplus\un\mpoint$$
Thus $C_2\times P\cong P\oplus P$, and the first assertion follows. The second one follows from Example~\ref{D8Q8}, which shows that in $\mathcal{R}_2$
$$D_8\cong \un\oplus4\cdot\partial C_2\cong 5\cdot \un\mpoint$$
Hence $D_8\times P\cong 5\cdot P$, thus $F(D_8\times P)\cong F(P)^{\oplus 5}$. 
\findemo
\begin{mth}{Proposition} \label{Q8 involution}The edge $\partial Q_8$ is an involution~: more precisely
$$\partial Q_8\times\partial Q_8=\partial C_2\cong\un\mpoint$$
\end{mth}
\pf Indeed $Q_8\diamond Q_8=C_2$, and $\nu_{Q_8,Q_8}=1$, by Theorem~\ref{tensor}.\findemo
\begin{rem}{Remark} The ``action" of this involution on the edges of the other Roquette 2-groups (that is, different from $\un$, $C_2$, and $Q_8$) is as follows~: it stabilizes cyclic and semidihedral groups, and exchanges dihedral and generalized quaternion groups. More precisely, it follows from Corollary~\ref{edge product} that
\begin{eqnarray*}
\forall n\geq 2,\;\partial Q_8\times \partial C_{2^n}&=&\partial C_{2^n}\\
\forall n\geq 4,\;\partial Q_8\times \partial D_{2^n}&=&\partial Q_{2^n}\\
\forall n\geq 4,\;\partial Q_8\times \partial Q_{2^n}&=&\partial D_{2^n}\\
\forall n\geq 4,\;\partial Q_8\times \partial SD_{2^n}&=&\partial SD_{2^n}\\
\end{eqnarray*}\vspace{-5ex}
\end{rem}
\npar By Theorem~\ref{Roquette category}, any finite $p$-group is isomorphic to a direct sum of edges of Roquette $p$-groups in the category $\mathcal{R}_p$. The following result shows that the summands of such an arbitrary direct sum are unique, up to group isomorphism, with the possible exception of the isomorphism $\un=\partial \un\cong \partial C_2$ of Proposition~\ref{edge C2}~:
\begin{mth}{Proposition} \label{unique}Let $\mathcal{S}$ and $\mathcal{T}$ be finite sequences of Roquette $p$-groups, such that there exists an isomorphism 
\begin{equation}\label{iso}
\dirsum{S\in\mathcal{S}}\partial S\cong \dirsum{T\in\mathcal{T}}\partial T
\end{equation}
in the category $\mathcal{R}_p$. If $p=2$, replace any occurrence of $C_2$ in $\mathcal{S}$ and $\mathcal{T}$ by the trivial group, which does not change the existence of the isomorphism~\ref{iso}, by Proposition~\ref{edge C2}.\par
Then there exists a bijection $\varphi:\mathcal{S}\to\mathcal{T}$ such that the groups $S$ and $\varphi(S)$ are isomorphic, for any $S\in\mathcal{S}$.
\end{mth}
\pf By~\cite{fonctrq}, the simple biset functors $S_{R,\F_p}$, where $R$ is a Roquette $p$-group different from $C_p$, are rational biset functors. Moreover, if $|R|\geq p^2$, then for any finite $p$-group $P$, the dimension of $S_{R,\F_p}(P)$ is equal to the number of groups $S$ in a genetic basis of $P$ such that $\sur{N}_P(S)\cong R$. On the other hand, the $\F_p$-dimension of $S_{\un,\F_p}(P)$ is equal to the number of groups $S$ in a genetic basis of $P$ such that $|\sur{N}_P(S)|\leq p$.\par
The functor $S_{R,\F_p}$ extends to an additive functor from $\mathcal{R}_p$ to the category of $\F_p$-vector spaces, and the value of this functor at the edge $\partial P$ is by definition equal to $\partial S_{R,\F_p}(P)$. If $|R|\geq p^2$, then the $\F_p$-dimension of $\partial S_{R,\F_p}(P)$ is equal to the number of groups $S$ in a genetic basis of $P$ such that $\sur{N}_P(S)\cong R$ and $S\cap Z(P)=\un$. In particular, if $P$ itself is a Roquette group, then 
$$\dim_{\F_p}S_{R,\F_p}(\partial P)=\dim_{\F_p}\partial S_{R,\F_p}(P)=\left\{\begin{array}{cl}1&\hbox{if }\;P\cong R\\0&\hbox{otherwise}\end{array}\right.$$
Applying the functor $S_{R,\F_p}$ to the isomorphism~\ref{iso}, this implies that the number of terms in the sequence $\mathcal{S}$ which are isomorphic to $R$ is equal to the corresponding number in the sequence $\mathcal{T}$.\par
Similarly, for any finite $p$-group $P$, the $\F_p$-dimension of $\partial S_{\un,\F_p}(P)$ is equal to the number of groups $S$ in a genetic basis of $P$ such that $\sur{N}_P(S)\leq p$ and $S\cap Z(P)=\un$. If $P$ itself is a Roquette group, this gives
$$\dim_{\F_p}S_{\un,\F_p}(\partial P)=\dim_{\F_p}\partial S_{\un,\F_p}(P)=\left\{\begin{array}{cl}1&\hbox{if}\;P\cong C_p\\1&\hbox{if }\;P\cong \un\\0&\hbox{otherwise}\end{array}\right.$$
Hence the number of terms in the sequence $\mathcal{S}$ which are isomorphic to $\un$ or $C_p$ is equal to the corresponding number in the sequence $\mathcal{T}$. If $p=2$, there are no $S$ in $\mathcal{S}\cup\mathcal{T}$ such that $S\cong C_2$, by assumption. It follows that for any Roquette $p$-group $R$, the number of terms in the sequence $\mathcal{S}$ which are isomorphic to $R$ is equal to the corresponding number in the sequence $\mathcal{T}$. The proposition follows in this case.\par
If $p>2$, the above argument shows that
$$\oplusb{S\in\mathcal{S}}{|S|\geq p^2}\partial S\cong\oplusb{T\in\mathcal{T}}{|T|\geq p^2}\partial T\mpoint$$
Let $M$ denote this direct sum. The isomorphism~\ref{iso} can be rewritten as
\begin{equation}\label{iso2}
m_\un\un\oplus m_{C_p}\partial C_p\oplus M\cong
n_\un\un\oplus n_{C_p}\partial C_p\oplus M\mvirg
\end{equation}
for some integers $m_\un, m_{C_p}, n_\un, n_{\C_p}$ such that $m_\un+m_{C_p}=n_\un+n_{C_p}$.\par
Now, let $\zeta$ be a primitive (i.e. non-trivial, since $p$ is prime) character $(\Z/p\Z)^\times\to \C$. Such a character exists since $p>2$. The functor $S_{C_p,\zeta}$ is a rational $p$-biset functor, as it is a summand of $\C R_\C$ (\cite{bisetfunctors} Corollary~7.3.5). Applying this functor to the isomorphism~\ref{iso2} and taking dimensions gives
$$m_{C_p}+\dim_{\C}S_{C_p,\zeta}(M)=n_{C_p}+\dim_{\C}S_{C_p,\zeta}(M)\mvirg$$
since $S_{C_p,\zeta}(\un)=0$ and $S_{C_p,\zeta}(C_p)=\partial S_{C_p,\zeta}(C_p)\cong \C$. \par
It follows that $m_{C_p}=n_{C_p}$, hence $m_\un=n_\un$, which completes the proof.\findemo
\vspace{2ex}
\pagebreak[3]
\begin{mth}{Corollary} \label{simplifiable}Let $X$, $Y$, and $Z$ be objects of $\mathcal{R}_p$, isomorphic to direct sums of edges of Roquette $p$-groups. 
\begin{enumerate}
\item If $X\oplus Z\cong Y\oplus Z$ in $\mathcal{R}_p$, then $X\cong Y$.
\item If $n$ is a positive integer, and if $n\cdot X\cong n\cdot Y$ in $\mathcal{R}_p$, then $X\cong Y$.
\end{enumerate}
\end{mth}
\pf Decompose $X$ as $X\cong\dirsum{R}n_R(X)\cdot\partial R$, where $R$ runs through the set of isomorphism classes of Roquette $p$-groups, and the function $R\mapsto n_R(X)\in\nolinebreak\N$ has finite support. Choose similar decompositions $Y\cong\dirsum{R}n_R(Y)\cdot\partial R$ and $Z\cong\dirsum{R}n_R(Z)\cdot\partial R$.\par
\pagebreak[3]
For Assertion~1, if $p>2$, it follows from Proposition~\ref{unique} that 
$$n_R(X)+n_R(Z)=n_R(Y)+n_R(Z)\mvirg$$
for each $R$. Thus $n_R(X)=n_R(Y)$ for each $R$, hence $X\cong Y$ in $\mathcal{R}_p$.\par
If $p=2$, and if $R$ is a Roquette $p$-group different from $\un$ and $C_2$, Proposition~\ref{unique} shows that $n_R(X)+n_R(Z)=n_R(Y)+n_R(Z)$, hence $n_R(X)=n_R(Y)$. Proposition~\ref{unique} also implies that
$$n_\un(X)+n_{C_2}(X)+n_\un(Z)+n_{C_2}(Z)=n_\un(Y)+n_{C_2}(Y)+n_\un(Z)+n_{C_2}(Z)\mvirg$$
whence $n_\un(X)+n_{C_2}(X)=n_\un(Y)+n_{C_2}(Y)$, and $X\cong Y$ in $\mathcal{R}_2$ again, since $\un\cong C_2$.\par
The proof of Assertion~2 is similar~: if $p>2$, Proposition~\ref{unique} shows that $n\, n_R(X)=n\, n_R(Y)$, for any $R$, thus $n_R(X)=n_R(Y)$, and $X\cong Y$. And if $p=2$, the conclusion $n_R(X)=n_R(Y)$ is valid for $R$ different from $\un$ and $C_2$. Moreover $n \big(n_\un(X)+n_{C_2}(X)\big)= n \big(n_\un(Y)+n_{C_2}(Y)\big)$, hence $n_\un(X)+n_{C_2}(X)= n_\un(Y)+n_{C_2}(Y)$, and $X\cong Y$, since $\un\cong C_2$.  \findemo
In the case of the decomposition of a $p$-group as a direct sum of edges of Roquette groups, the above isomorphism $\partial C_2\cong \partial\un$ doesn't matter, and the decomposition is unique~:
\begin{mth}{Proposition} \label{isomorphism}Let $P$ and $Q$ be finite $p$-groups. The following assertions are equivalent~:
\begin{enumerate}
\item The groups $P$ and $Q$ are isomorphic in the category $\mathcal{R}_p$.
\item There exist genetic bases  $\mathcal{B}_P$ and $\mathcal{B}_Q$ of $P$ and $Q$, respectively, and a bijection $\sigma:\mathcal{B}_P\stackrel{\cong}{\longrightarrow}\mathcal{B}_Q$ such that
\begin{equation}\label{isobases}
\forall S\in\mathcal{B}_P,\;\;N_{Q}\big(\sigma(S)\big)/\sigma(S)\cong N_P(S)/S\mpoint
\end{equation}
\item For any genetic bases $\mathcal{B}_P$ and $\mathcal{B}_Q$ of $P$ and $Q$, respectively, there exists a bijection $\sigma:\mathcal{B}_P\stackrel{\cong}{\longrightarrow}\mathcal{B}_Q$ such that~\ref{isobases} holds.
\end{enumerate}
\end{mth}
\pf 
Assertion~2 implies Assertion~1 by Theorem~\ref{Roquette category}. Now suppose that Assertion~1 holds. Then in particular $F(P)\cong F(Q)$, for any rational $p$-biset functor $F$. Let $\mathcal{B}_P$ and $\mathcal{B}_Q$ be genetic bases of $P$ and $Q$, respectively. If $R$ is a Roquette $p$-group, set
$$m_P(R)=\big|\{S\in\mathcal{B}\mid N_P(S)/S\cong R\}\big|\mvirg$$
and define similarly $m_{Q}(R)$ for the group $Q$. The integers $m_P(R)$ and $m_Q(R)$ do not depend on the choices of the genetic bases $\mathcal{B}_P$ and $\mathcal{B}_Q$.\par
If $R$ is not isomorphic to $C_p$, then the simple functor $S_{R,\F_p}$ is rational. Moreover, the $\F_p$-dimension of $S_{R,\F_p}(P)$ is equal to $m_P(R)$ if $|R|>p$, and to $1+m_{P}(C_p)$, if $R=\un$. Since $P$ is the only element $S$ of $\mathcal{B}$ such that $N_P(S)/S=\un$, it follows $m_P(\un)=1$, and then $m_P(R)=m_{Q}(R)$ for any Roquette $p$-group $R$. Assertion~2 follows. The equivalence of Assertions 2 and 3 follows from Theorem~\ref{genetic recall}.\findemo
\begin{rem}{Examples} \label{extra}\mpn
$\bullet$ Let $p>2$, and let $X^+$ (resp. $X^-$) denote the extraspecial $p$-group of order $p^3$ and exponent $p$ (resp. $p^2$). Then $X^+\cong X^-$ in~$\mathcal{R}_p$, for if~$P$ is one of these groups, each genetic basis of $P$ consists of $S=P$, for which $N_P(S)/S=\un$, of the $p+1$ subgroups $S$ of index $p$ in $P$, for which $N_P(S)=P/S\cong C_p$, and an additional non-normal genetic subgroup~$S$ such that $N_P(S)/S\cong C_p$. In other words
$$X^+\cong X^-\cong\un\oplus(p+2)\cdot\partial C_p$$
in the category $\mathcal{R}_p$.  \mpn
$\bullet$ Similar examples exist for $p=2$~: if $P$ is one of the groups labelled 6 or~7 in the GAP list of groups of order 32 (see \cite{GAP4}), with respective structure $((C_4\times C_2)\rtimes C_2)\rtimes C_2$ and $(C_8\rtimes C_2)\rtimes C_2$, then in any genetic basis of $P$, there is a unique group $S(=P)$ such that $N_P(S)/S=\un$, there are 6 groups $S$ such that $N_P(S)/S\cong C_2$, and 2 groups $S$ such that $N_P(S)/S\cong C_4$.\mpn
$\bullet$ \label{different orders}Some 2-groups {\em with different orders} may become isomorphic in the category~$\mathcal{R}_2$~: using GAP, one can show that the elementary abelian group of order~16 is isomorphic to each of the groups labelled 134, 138, and 177 in GAP's list of groups of order 64. These groups have respective structure 
$$\big((C_4\times C_4)\rtimes C_2)\rtimes C_2,\;\;\Big(\big((C_4\times C_2)\rtimes C_2\big)\rtimes C_2\Big)\rtimes C_2,\;\;\hbox{and}\;\;(C_2\times D_{16})\rtimes C_2\mpoint$$
\vspace{-3ex}
\end{rem}
I couldn't find any similar example for $p>2$. In this case however, the following result characterizes those $p$-groups which become isomorphic in the category $\mathcal{R}_p$~:
\begin{mth}{Proposition} \label{isomorphism Roquette} Let $p$ be a prime number, and let $P$ and $Q$ be finite $p$-groups. \begin{enumerate}
\item If $P\cong Q$ in the category $\mathcal{R}_p$, then the $\Q$-algebras $Z\Q P$ and $Z\Q Q$ are isomorphic.
\item If $p>2$, and if $Z\Q P$ and $Z\Q Q$ are isomorphic $\Q$-algebras, then $P\cong Q$ in the category $\mathcal{R}_p$.
\end{enumerate}
\end{mth}
\pf Let $\mathcal{G}$ be a genetic basis of $P$. For $S\in\mathcal{G}$, let $V(S)$ denote the corresponding simple $\Q P$-module, defined by
$$V(S)=\Indinf_{\sur{N}_P(S)}^P\Phi_{\sur{N}_P(S)}\mpoint$$
The multiplicity $v_S$ of $V(S)$ in the $\Q P$-module $\Q P$ is equal to
$$v_S=\frac{\dim_\Q V(S)}{\dim_\Q\End_{\Q P}\big(V(S)\big)}\mpoint$$
As $S$ is a genetic subgroup of $P$, there is an isomorphism of (skew-)fields
$$\End_{\Q P}\big(V(S)\big)\cong \End_{\Q \sur{N}_P(S)}\big(\Phi_{\sur{N}_P(S)}\big)\mpoint$$
It follows that there is an isomorphism of $\Q$-algebras
$$\Q P\cong\prod_{S\in\mathcal{G}}M_{v_S}\Big(\End_{\Q \sur{N}_P(S)}\big(\Phi_{\sur{N}_P(S)}\big)\Big)\mpoint$$
Hence
$$Z\Q P\cong\prod_{S\in\mathcal{G}}Z\Big(\End_{\Q \sur{N}_P(S)}\big(\Phi_{\sur{N}_P(S)}\big)\Big)\mpoint$$
It shows that the isomorphism type of the $\Q$-algebra $Z\Q P$ depends only on the genetic basis $\mathcal{G}$~: more precisely, it is determined by the isomorphism type of $P$ in $\mathcal{R}_p$. This proves Assertion~1.\par
Now if $p>2$, the group $\sur{N}_P(S)$ is cyclic, of order $p^{m_S}$, say. By Example~\ref{endo Phi}, there is an isomorphism of (skew-)fields
$$\End_{\Q \sur{N}_P(S)}\big(\Phi_{\sur{N}_P(S)}\big)\cong \Q (\zeta_{p^{m_S}})\mvirg$$
where $\zeta_{p^{m_S}}$ is a primitive root of unity of order $p^{m_S}$. Hence
$$Z\Q P\cong\prod_{S\in\mathcal{G}}\Q(\zeta_{p^{m_S}})\mpoint$$
Similarly, if $\mathcal{H}$ is a genetic basis of $Q$
$$Z\Q Q\cong\prod_{T\in\mathcal{H}}\Q(\zeta_{p^{n_T}})\mvirg$$
where $p^{n_T}=|\sur{N}_Q(T)|$. \par
Let $l$ be an integer bigger than all the $m_S$'s, for $S\in\mathcal{G}$, and all the $n_T$'s, for $T\in\mathcal{H}$. Set $K=\Q(\zeta_{p^l})$, and let $G$ be the Galois group of $K$ over $\Q$. By Galois theory (\cite{szamuely} Theorem 1.5.4 and Remark 1.5.5), the $\Q$-algebras $Z\Q P$ and $Z\Q Q$ are isomorphic if and only if there is an isomorphism of $G$-sets
$$\Hom_{\rm alg}(Z\Q P,K)\cong \Hom_{\rm alg}(Z\Q Q,K)\mpoint$$ 
When $r\leq l$ is an integer, let $G_r$ denote the Galois group of $K$ over $\Q(\zeta_{p^{r}})$. Then the $G$-set $\Hom_{\rm alg}(Z\Q P,K)$ is isomorphic to
$$\mathop{\bigsqcup}_{S\in\mathcal{G}}\limits G/G_{n_S}\mpoint$$
The isomorphism $Z\Q P\cong Z\Q Q$ implies that for any $r\leq l$, the number of $S\in\mathcal{G}$ such that $\sur{N}_P(S)$ has order~$p^r$ is equal to the number of $T\in \mathcal{H}$ such that $\sur{N}_Q(T)$ has order $p^r$. Now Assertion~2 of the proposition follows from Proposition~\ref{isomorphism}.\findemo
\begin{rem}{Remark} \label{not true for p=2}Assertion~2 of Proposition~\ref{isomorphism Roquette} is not true for $p=2$~: let $P=D_8$ and $Q=Q_8$ denote a dihedral group of order 8 and a quaternion group of order 8, respectively. By Example~\ref{D8Q8}, in a genetic basis of $P$, there is one group $S$ such that $N_P(S)/S=\un$ (namely $S=P$), and 4 subgroups $S$ such that $N_P(S)/S\cong C_2$ (the 3 subgroups of index 2 in $P$, and a non-central subgroup of order 2 of $P$). It follows easily that
$$\Q D_8\cong \Q\oplus\Q\oplus\Q\oplus\Q\oplus M_2(\Q)\mpoint$$
On the other hand, a genetic basis of $Q$ contains one subgroup $S$ such that $N_Q(S)/S=\un$ (namely $S=Q$), 3 subgroups $S$ such that $N_Q(S)/S\cong C_2$ (the 3 subgroups of index 2 in $Q$), and one subgroup $S$ such that $N_Q(S)/S\cong Q_8$ (the trivial subgroup of $Q$). Hence
$$\Q Q_8\cong \Q\oplus\Q\oplus\Q\oplus\Q\oplus \mathbb{H}_\Q\mvirg$$
where $\mathbb{H}_\Q$ is the field of quaternions over $\Q$. Then
$$Z \Q D_8\cong \Q^5\cong Z\Q Q_8\mpoint$$
But $D_8$ and $Q_8$ are not isomorphic in $\mathcal{R}_2$, by Proposition~\ref{isomorphism}.
\end{rem}
\masubsect{Genetic bases of direct products} Theorem~\ref{tensor} yields a way to compute a genetic basis of a direct products of $p$-groups. More precisely:
\begin{mth}{Theorem} \label{genetic basis product}Let $P$ and $Q$ be finite $p$-groups, let $\mathcal{B}_P$ be a genetic basis of $P$, and let $\mathcal{B}_Q$ be a genetic basis of $Q$. 
\begin{enumerate}
\item For each pair $(S,T)\in\mathcal{B}_P\times\mathcal{B}_Q$, let $\sur{R}$ be a centrally diagonal genetic subgroup of $\sur{N}_P(S)\times \sur{N}_Q(T)$, and let
$$R=\big\{(x,y)\in N_P(S)\times N_Q(T)\mid (xS,yT)\in\sur{R}\big\}\mpoint$$
\nopagebreak
Then $R$ is a genetic subgroup of $P\times Q$, such that 
$$\sur{N}_{P\times Q}(R)\cong \sur{N}_P(S)\diamond \sur{N}_Q(T)\mpoint$$
\pagebreak[3]
\item For $(S,T)\in\mathcal{B}_P\times\mathcal{B}_Q$, let $\mathcal{E}_{S,T}$ denote the set of subgroups $R$ obtained in Assertion~1, when $\sur{R}$ runs through a set of representatives of centrally diagonal genetic subgroups of $\sur{N}_P(S)\times \sur{N}_Q(T)$, for the relation $\bizlie{\sur{N}_P(S)\times \sur{N}_Q(T)}$, as described in Theorem~\ref{tensor}. \par
Then the sets $\mathcal{E}_{S,T}$ consist of mutually inequivalent genetic subgroups of $P\times Q$, for the relation $\bizlie{P\times Q}$, and the (disjoint) union
$$\mathcal{B}_{P\times Q}=\mathop{\bigsqcup}_{(S,T)\in\mathcal{B}_P\times\mathcal{B}_Q}\mathcal{E}_{S,T}$$
is a genetic basis of $(P\times Q)$.
\end{enumerate}
\end{mth}
\pf Assertion~1 is straightforward if the group $\sur{N}_P(S)$ is trivial, i.e. if $S=P$, or if the group $\sur{N}_Q(T)$ is trivial, i.e. if $T=Q$. So we can assume that $S<P$ and $T<Q$. \par
By Theorem~\ref{genetic product Roquette}, the group $\sur{R}$ is diagonal in $\sur{N}_P(S)\times\sur{N}_Q(T)$. It follows that $k_1(R)=S$, $k_2(R)=T$, and $\sur{R}=R/(S\times T)$. This implies in particular that $N_{P\times Q}(R)\leq N_P(S)\times N_Q(T)$. More precisely
$$N_{P\times Q}(R)=\big\{(a,b)\in N_P(S)\times N_Q(T)\mid (aS,bT)\in N_{\sur{N}_P(S)\times\sur{N}_Q(T)}(\sur{R})\big\}\mvirg$$
and the map $(a,b)\mapsto (aS,bT)$ induces a group isomorphism
$$N_{P\times Q}(R)/R\cong N_{\sur{N}_P(S)\times\sur{N}_Q(T)}(\sur{R})/\sur{R}\mpoint$$
It follows that $N_{P\times Q}(R)/R$ is a Roquette group, and by Theorem~\ref{genetic product Roquette} again, and Notation~\ref{diamond}
$$N_{P\times Q}(R)/R\cong\sur{N}_P(S)\diamond \sur{N}_Q(T)\mpoint$$
Let $\widehat{S}\geq S$ denote the subgroup of $N_P(S)$ such that $\widehat{S}/S$ is the unique central subgroup of order $p$ of the Roquette group $\sur{N}_P(S)$. Define $\widehat{T}\geq T$ similarly, and let $\widehat{R}/R$ be the unique central subgroup of order $p$ of $N_{P\times Q}(R)/R$. Then
$$\widehat{R}=(\widehat{S}\times\un)R=(\un\times\widehat{T})R\mpoint$$
Let $(x,y)\in P\times Q$ such that $R^{(x,y)}\cap \widehat{R}\leq R$. Intersecting this inclusion with $P\times \un$ gives
$$(S^x\cap\widehat{S})\times \un\leq S\times \un\mvirg$$
thus $S^x\cap \widehat{S}\leq S$, and it follows that $x\in N_P(S)$, since $S$ is an expansive subgroup of $P$. Similarly, intersecting the inclusion $R^{(x,y)}\cap \widehat{R}\leq R$ with $\un\times Q$ gives $T^y\cap\widehat{T}\leq T$, hence $y\in N_Q(T)$.\par
Now $S\times T\leq R^{(x,y)}\cap \widehat{R}\leq R$, and taking the quotient by $S\times T$ gives
$$\sur{R}^{(xS,yT)}\cap \big(\widehat{R}/(S\times T)\big)\leq \sur{R}\mpoint$$
As $\sur{R}$ is a genetic subgroup of $\sur{N}_P(S)\times\sur{N}_Q(T)$, it follows that $\sur{R}^{(xS,yT)}$ is equal to $\sur{R}$, hence $R^{(x,y)}=R$. Thus $R$ is an expansive subgroup of $P\times Q$. Since $N_{P\times Q}(R)/R$ is a Roquette group, the group $R$ is a genetic subgroup of $P\times Q$, and this completes the proof of Assertion~1.\par
For Assertion~2, let $(S,T)$ and $(S',T')$ in $\mathcal{B}_P\times\mathcal{B}_Q$, and let $R\in\mathcal{E}_{S,T}$ and $R'\in\mathcal{E}_{S',T'}$, such that $R\bizlie{P\times Q}R'$. It means that there exists $(x,y)\in P\times Q$ such that
\begin{equation}\label{intersect}
R^{(x,y)}\cap \widehat{R}'\leq R'\mvirg\;\;^{(x,y)}R'\cap\widehat{R}\leq R\mpoint
\end{equation}
Intersecting these two inclusions with $P\times \un$ gives
$$S^x\cap \widehat{S}'\leq S'\mvirg\;\;^xS'\cap \widehat{S}\leq S\mpoint$$
Hence $S'\bizlie{P}S$, thus $S'=S$, since $S$ and $S'$ are in the same genetic basis of $P$. Moreover $x\in N_P(S)$. Similarly, intersecting~\ref{intersect} with $\un\times Q$ implies $T=T'$, and $y\in N_Q(T)$. \par
Quotienting the inclusions~\ref{intersect} by $(S\times T)$ gives that $\sur{R}'\bizlie{\sur{N}_P(S)\times\sur{N}_Q(T)} \sur{R}$. Hence $R'=R$, as was to be shown.\par
Now setting
$$\mathcal{B}_{P\times Q}=\mathop{\bigsqcup}_{(S,T)\in\mathcal{B}_P\times\mathcal{B}_Q}\mathcal{E}_{S,T}$$
yields a set of genetic subgroups of $P\times Q$, which are inequivalent to one other for the relation $\bizlie{P\times Q}$. But
$$|\mathcal{E}_{S,T}|=\nu_{\sur{N}_P(S),\sur{N}_Q(T)}\mvirg$$
and $N_{P\times Q}(R)/R\cong \sur{N}_P(S)\diamond\sur{N}_Q(T)$, for any $R\in\mathcal{E}_{S,T}$
it follows that
\begin{eqnarray*}
\dirsum{R\in\mathcal{B}_{P\times Q}}\partial N_{P\times Q}(R)/R&\cong&\oplusb{S\in\mathcal{B}_P}{T\in\mathcal{B}_Q}\nu_{\sur{N}_P(S),\sur{N}_Q(T)}\partial\big(\sur{N}_P(S)\diamond\sur{N}_Q(T)\big)\\
&\cong&\big(\dirsum{S\in\mathcal{B}_P}\partial \sur{N}_P(S)\big)\times\big(\dirsum{T\in\mathcal{B}_Q}\partial \sur{N}_Q(T)\big)\\
&\cong& P\times Q\mpoint
\end{eqnarray*}
In particular, the rank $l_\Q(P\times Q)$ of the group $R_\Q(P\times Q)$ is equal to $|\mathcal{B}_{P\times Q}|$. Since $\mathcal{B}_{P\times Q}$ is contained in a genetic basis of $P\times Q$, which has cardinality $l_\Q(P\times Q)$, it follows that $\mathcal{B}_{P\times Q}$ is a genetic basis of $P\times Q$.\findemo
\begin{rem}{Remark} Theorem~\ref{genetic basis product} does not mean that any genetic subgroup of $P\times Q$ can be obtained by the construction of Assertion~1. For example, if $[P,P]\leq Z(P)$ and $Z(P)$ is cyclic, then the diagonal $R=\Delta(P)$ is a genetic subgroup of $P\times P$, by Example~\ref{diagonal}. But $k_1(R)=\un$ is not a genetic subgroup of $P$, if $P$ is not a Roquette group.
\end{rem}
\masubsect{Example of application} As explained in Example~\ref{D8Q8}, the dihedral group $D_8$ splits as
\begin{equation}\label{iso D_8}D_8\cong \un\oplus 4\partial C_2
\end{equation}
in the category $\mathcal{R}_2$. By Proposition~\ref{edge C2}, it follows that $D_8\cong 5\cdot\un$ in $\mathcal{R}_2$. Hence $(D_8)^n\cong 5^n\cdot\un$, for any $n\in\N$. In particular, if $F$ is a rational $2$-biset functor such that $F(\un)=\zero$, then $F\big((D_8)^n\big)=\zero$. Hence $F(P)=\zero$ for any quotient of a direct product of copies of $D_8$, by Remark~\ref{factor group}.\par
Actually, one can be more precise~: since $\partial C_2\times \partial C_2\cong\partial C_2$ by Corollary~\ref{edge product}, it follows that for any $n\in\N$
\begin{equation}\label{D2n}(D_8)^n\cong\dirsum{i=0}^{n}\binom{n}{i}4^{i}\cdot(\partial C_2)^i=\un\oplus \dirsum{i=1}^{n}\binom{n}{i}4^{i}\cdot\partial C_2\cong \un\oplus (5^n-1)\cdot\partial C_2\mpoint
\end{equation}
It means that a genetic basis of the group $P=(D_8)^n$ is made of the group $S=P$, for which $N_P(S)/S=\un$, and of $(5^n-1)$ subgroups $S$ for which $N_P(S)/S\cong C_2$. \par
In particular, by Theorem~9.5 of~\cite{dadegroup} (or Corollary~12.10.3 of~\cite{bisetfunctors}), the Dade group of $P$ is torsion free, and so is the Dade group of any factor group of~$P$, by Remark~\ref{factor group} again. It shows that the Dade group of a central product of any number of copies of~$D_8$ is torsion free (see Theorem~\ref{product of dihedral groups} for a generalization of this result)~: this was proved by Nadia Mazza and me (Theorem~9.2 of~\cite{boma}). However, the above argument cannot be considered as a new proof of this result, since Theorem~9.5 of~\cite{dadegroup} relies on Theorem~9.2 of~\cite{boma}.
\masubsect{Edges of central products}
Let $P$ and $Q$ be non-trivial finite $p$-groups. Recall that {\em a central product} $P*_\varphi Q$ of $P$ and $Q$ is by definition a group of the form $(P\times Q)/\rDelta_\varphi(Z_P)$, where $Z_P$ is a central subgroup of order~$p$ of $P$, and $\varphi : Z_P\hookrightarrow Z(Q)$ is some isomorphism from $Z_P$ to some central subgroup $Z_Q$ of $Q$. \par
In the case where $p=2$ and the groups $P$ and $Q$ both have cyclic center, the group $Z_P$ is unique, as well as the morphism $\varphi$, so the central product is simply denoted by $P*Q$ in this case.
\begin{mth}{Proposition} \label{edge of central products}Let $p$ be a prime number, and let $P$ and $Q$ be non-trivial finite $p$-groups. Let $Z_P$ (resp. $Z_Q$) denote a central subgroup of order~$p$ of $P$ (resp. $Q$).
\begin{enumerate}
\item If one of the groups $Z(P)$ or $Z(Q)$ is non-cyclic, or if $|Z(P)|>p$  and $|Z(Q)|>p$, then $\partial (P*_\varphi Q)=0$ in $\mathcal{R}_p$, for any group isomorphism $\varphi:Z_P\to Z_Q$.
\item If $Z(P)$ and $Z(Q)$ are cyclic, and if moreover $Z(P)$  or $Z(Q)$ has order~$p$, then, in $\mathcal{R}_p$,
$$\dirsum{\varphi:Z_P\stackrel{\cong}{\to}Z_Q}\partial(P*_\varphi Q)\cong \partial P\times\partial Q\mpoint\vspace{-4ex}$$
\end{enumerate}
\end{mth}
\pf The center of the group $P*_\varphi Q$ is equal to $Z(P)*_\varphi Z(Q)$. It is cyclic if and only if both $Z(P)$ and $Z(Q)$ are cyclic, and if one of them has order~$p$. This proves Assertion~1.\par
For Assertion~2, suppose that $Z(P)$ and $Z(Q)$ are cyclic, and that one of them has order $p$. Then the subgroups $Z_P$ and $Z_Q$ are uniquely determined, and there are $p-1$ group isomorphisms $\varphi:Z_P\to Z_Q$. For each of them, the only central subgroup $Z_\varphi$ of order $p$ of $P*_\varphi Q$ is equal to $(Z_P\times Z_Q)/\rDelta_\varphi(Z_P)$, and 
$$(P*_\varphi Q)/Z_\varphi\cong (P\times Q)/(Z_P\times Z_Q)\cong\sur{P}\times\sur{Q}\mvirg$$
where $\sur{P}=P/Z_P$ and $\sur{Q}=Q/Z_Q$.\par
By Proposition~\ref{sum of edges of quotients}
\newlength{\longl}
\settowidth{\longl}{$\dirsum{(Z_P\times Z_Q)\leq N\normal (P\times Q)}$}
\begin{eqnarray*}
P\times Q&\cong&\makebox[\longl]{$\dirsum{\un\leq N\normal (P\times Q)}$}\partial\big((P\times Q)/N\big)\\
P*_\varphi Q&\cong&\makebox[\longl]{$\dirsum{\Delta_\varphi(Z_P)\leq N\normal (P\times Q)}$}\partial\big((P\times Q)/N\big)\\
\sur{P}\times Q&\cong&\makebox[\longl]{$\dirsum{(Z_P\times\un)\leq N\normal (P\times Q)}$}\partial\big((P\times Q)/N\big)\\
P\times \sur{Q}&\cong&\makebox[\longl]{$\dirsum{(\un\times Z_Q)\leq N\normal (P\times Q)}$}\partial\big((P\times Q)/N\big)\\
\sur{P}\times \sur{Q}&\cong&\makebox[\longl]{$\dirsum{(Z_P\times Z_Q)\leq N\normal (P\times Q)}$}\partial\big((P\times Q)/N\big)\mpoint\\
\vspace{-2ex}\end{eqnarray*}
Set \vspace{-2ex}
\begin{equation}\label{S egale}S=\partial(P\times Q)\oplus\left(\dirsum{\varphi:Z_P\stackrel{\cong}{\to}Z_Q}(P*_\varphi Q)\right)\oplus(\sur{P}\times Q)\oplus(P\times \sur{Q})\mpoint
\end{equation}
\pagebreak[3]
Then $S$ is equal to the direct sum of the edges $\partial\big((P\times Q)/N\big)$, for $N\normal (P\times Q)$, with multiplicity $(p+1)$ if $N\geq (Z_P\times Z_Q)$, and multiplicity 1 otherwise. Hence
\begin{equation}\label{S egale2}S\cong (P\times Q)+p\cdot (\sur{P}\times \sur{Q}\big)\mpoint
\end{equation}
Now $\partial (P\times Q)=0$ by Corollary~\ref{zero edge}. Moreover, by Corollary~\ref{edge cyclic center}, for each $\varphi:Z_P\stackrel{\cong}{\to}Z_Q$,
$$P*_\varphi Q=\partial(P*_\varphi Q)\oplus (\sur{P}\times\sur{Q})\mpoint$$
But $P\cong\partial P\oplus \sur{P}$ and $Q\cong\partial Q\oplus \sur{Q}$, by Corollary~\ref{edge cyclic center} again. Replacing $P$, $Q$, and $P*_\varphi Q$ by these values in Equation~\ref{S egale} gives
\begin{equation}\label{S1}S\cong\left(\dirsum{\varphi:Z_P\stackrel{\cong}{\to}Z_Q}\partial(P*_\varphi Q)\right)\oplus (p+1)\cdot(\sur{P}\times\sur{Q})\oplus \big(\sur{P}\times(\partial Q)\big)\oplus \big((\partial P)\times \sur{Q}\big)\mpoint
\end{equation}
But replacing $P$ by $\partial P\oplus \sur{P}$ and $Q$ by $\partial Q\oplus \sur{Q}$ in Equation~\ref{S egale2} gives
\begin{equation}\label{S2}S\cong(\partial P\times\partial Q)\oplus (p+1)\cdot(\sur{P}\times\sur{Q})\oplus \big(\sur{P}\times(\partial Q)\big)\oplus \big((\partial P)\times \sur{Q}\big)\mpoint
\end{equation}
Comparing Equations~\ref{S1} and~\ref{S2} gives 
$$\dirsum{\varphi:Z_P\stackrel{\cong}{\to}Z_Q}\partial(P*_\varphi Q)\cong \partial P\times\partial Q\mvirg$$
by Corollary~\ref{simplifiable}. This completes the proof.\findemo
\begin{mth}{Corollary} \label{edge of 2-groups}\begin{enumerate}
\item Let $P$ and $Q$ be non-trivial finite 2-groups with cyclic center, and assume that $Z(P)$ or $Z(Q)$ has order 2. Then
$$\partial(P*Q)\cong \partial P\times \partial Q$$
in the category $\mathcal{R}_2$.
\item For each $i\in\{1,\ldots,n\}$, let $P_i$  be a finite 2-group with center $Z_i$ of order 2. Let $\rule{0ex}{4ex}\monast{i=1}{n}P_i=P_1*P_2*\cdots*P_n$ denote the central product of the groups $P_i$. Then
$$\monast{i=1}{n}P_i\cong \prod_{i=1}^n(\partial P_i)\oplus \prod_{i=1}^n(P_i/Z_i)$$
in the category $\mathcal{R}_2$.
\item In particular, for any positive integer $n$, and any integer $m\geq 4$, there are isomorphisms
\begin{eqnarray*}
(D_{2^m})^{*n}&\cong&2^{(n-1)(m-3)}\cdot\partial D_{2^m}\oplus (D_{2^{m-1}})^n\\
(SD_{2^m})^{*n}&\cong&2^{(n-1)(m-3)}\cdot\partial SD_{2^m}\oplus (D_{2^{m-1}})^n\\
(Q_{2^m})^{*n}&\cong&\left\{\begin{array}{l}2^{(n-1)(m-3)}\cdot\partial D_{2^m}\oplus (D_{2^{m-1}})^n\;\hbox{if $n$ is even}\\
2^{(n-1)(m-3)}\cdot\partial Q_{2^m}\oplus (D_{2^{m-1}})^n\;\hbox{if $n$ is odd}
\\\end{array}\right.
\end{eqnarray*}
in the category $\mathcal{R}_2$, where $P^{*n}$ denote the central product of $n$ copies of~$P$.
\end{enumerate}
\end{mth}
\pf For Assertion~1, the assumptions imply that there is a unique isomorphism $\varphi:Z_P\to Z_Q$. Hence there is only one term in the summation of Proposition~\ref{edge of central products}.\par
Assertion~2 follows from Corollary~\ref{edge cyclic center}, which gives an isomorphism
$$\monast{i=1}{n}P_i\cong \partial \Big(\monast{i=1}{n}P_i\Big)\oplus \Big(\Big(\monast{i=1}{n}P_i\Big)/Y\Big)\mvirg$$
where $Y$ is the unique central subgroup of order $p$ in $\monast{i=1}{n}P_i$. Now an easy induction argument, using Assertion~1, shows that $\partial \Big(\monast{i=1}{n}P_i\Big)\cong\prod_{i=1}^n\limits(\partial P_i)$, and that $\Big(\monast{i=1}{n}P_i\Big)/Y\cong \prod_{i=1}^n\limits(P_i/Z_i)$.\par
Finally, when $P$ is one of the groups $D_{2^m}$, $SD_{2^m}$, or $SD_{2^m}$, then $P/Z\cong D_{2^{m-1}}$. Now Assertion~3 follows from Assertion~2 and from an easy induction argument using Corollary~\ref{edge product}.  \findemo
\begin{rem}{Remark} It follows from Assertion~3 that when $n$ is even, the groups $(D_{2^m})^{*n}$ and $(Q_{2^m})^{*n}$ are isomorphic in the category $\mathcal{R}_2$~: it is actually easy to check that they are isomorphic as {\em groups}.
\end{rem}
\begin{rem}{Example} \label{edges}From Corollary~\ref{edge product} and Assertion~1, it follows that~:
\begin{itemize} 
\item $\partial\big((D_8)^{*n}\big)\cong\partial C_2$.
\item $\partial\big((Q_8)^{*n}\big)\cong\left\{\begin{array}{cl}\partial C_2&\hbox{if $n$ is even}\\\partial Q_8&\hbox{if $n$ is odd}\end{array}\right.$
\item $\partial\big((SD_{2^m})^{*n}\big)\cong 2^{(n-1)(m-3)}\cdot\partial SD_{2^m}$, for $m\geq 4$.
\end{itemize}
More generally, if $P$ is any central product of groups isomorphic to $D_8$ or $Q_8$, that is, if $P$ is an extraspecial 2-group, then $\partial P\cong \partial C_2$, or $\partial P\cong\partial Q_8$. In particular (see Example~\ref{representation edge}), we recover the well known fact that $P$ has a unique faithful rational irreducible representation. But there is more~: let $Q$ be a non trivial 2-group. If the center of $Q$ is not cyclic, then the center of any central product $P*Q$ is not cyclic, hence $\partial Q=\partial (P*Q)=0$ in $\mathcal{R}_2$. If the center of $Q$ is cyclic, then there is a unique central product $P*Q$. By Theorem~\ref{Roquette category}, there is a finite sequence $\mathcal{S}$ of Roquette 2-groups such that
$$\partial Q\cong\dirsum{R\in\mathcal{S}}\partial R$$
in the category $\mathcal{R}_2$. By Assertion~1 of Corollary~\ref{edge of 2-groups}, it follows that
$$\partial(P*Q)\cong \dirsum{R\in\mathcal{S}}(\partial P\times\partial R)\mpoint$$
Now $\partial P\cong \partial C_2$ or $\partial P\cong\partial Q_8$. In both cases, by Proposition~\ref{edge C2} and Proposition~\ref{Q8 involution}, the multiplication by $\partial P$ is a permutation of the edges of the Roquette $2$-groups. It follows that there is a sequence $\mathcal{S}'$ of Roquette 2 groups, of the same length as $\mathcal{S}$, such that
$$\partial (P*Q)\cong\dirsum{R\in\mathcal{S}'}\partial R\mpoint$$
In particular, for any field $K$ of characteristic 0, the groups $R_K(\partial P)=\partial R_K(P)$ and $R_K\big(\partial(P*Q)\big)=\partial R_K(P*Q)$ are free of the same rank, equal to the length of $\mathcal{S}$ or $\mathcal{S}'$.\par
Hence in any case, the groups $Q$ and $P*Q$ have the same number (possibly~0 if the center of $Q$ is not cyclic) of faithful irreducible representations over $K$, up to isomorphism.\mmp
Similarly, the last example above means in particular  that the group $(SD_{2^m})^{*n}$ admits $2^{(n-1)(m-3)}$ non isomorphic faithful rational irreducible representations.\par 
\end{rem}
\begin{rem}{Example} \label{extraspecial} Let $p$ be an odd prime, and let $P=X^\epsilon$ (where $\epsilon\in\{\pm 1\}$) be one of the extraspecial groups of order $p^3$ considered in Examples~\ref{extra}. Then $\partial P\cong\partial C_p$. Moreover $Z(P)$ has order $p$, and any automorphism of $Z(P)=Z_P$ can be extended to an automorphism of $P$. It follows that $P*_\varphi Q$ is independent (up to a group isomorphism) of the choice of an embedding $\varphi:P\hookrightarrow Z(Q)$, for any non-trivial $p$-group $Q$ with cyclic center, so we can denote this group by $P*Q$. \par
Now if $Q$ is a non-trivial $p$-group, and $\mathcal{B}$ is a genetic basis of $Q$, it follows from Theorem~\ref{Roquette category} that
$$\partial Q\cong\oplusb{S\in\mathcal{B}}{S\cap Z(Q)=\un}\partial \sur{N}_Q(S)\mvirg$$
and the right hand side is a direct sum of {\em non-trivial} Roquette $p$-groups. By Equation~\ref{product cyclic edges}, for any non trivial Roquette $p$-group $R$
$$\partial C_p\times\partial R\cong (p-1)\partial R\mpoint$$
Hence for any non-trivial $p$-group $Q$
$$\partial C_p\times\partial Q\cong (p-1)\partial Q\mpoint$$
It follows that if $Q$ is a non-trivial $p$-group with cyclic center, and if $P\cong X^\epsilon$, then
$$\dirsum{\varphi:Z_P\stackrel{\cong}{\to}Z_Q}\partial(P*_\varphi Q)\cong(p-1)\partial(P*Q)\cong \partial C_p\times\partial Q\cong (p-1)\partial Q\mvirg$$
hence $\partial(X^\epsilon*Q)\cong\partial Q$, by Corollary~\ref{simplifiable}. Note that this is also true (for any central product of $X^\epsilon$ with $Q$) if the center of $Q$ is non-trivial, since in this case the center of $X^\epsilon*Q$ is also non trivial, and then $\partial(X^\epsilon*Q)\cong\partial Q\cong 0$ in $\mathcal{R}_p$, by Corollary~\ref{zero edge}. \par
It follows easily by induction that if $P$ is any central product of groups isomorphic to $X^+$ or $X^-$, i.e. if $P$ is an extraspecial $p$-group, then $\partial P\cong \partial C_p$. We recover this way the well known fact that for $p>2$ also, extraspecial $p$-groups have a unique faithful rational irreducible representation. The same argument shows more generally that $\partial(P*Q)\cong\partial Q$ in $\mathcal{R}_p$, for any non-trivial $p$-group $Q$. In particular $Q$ and $P*Q$ have the same number of faithful irreducible representations over a given field $K$ of characteristic 0.
\end{rem}
\begin{mth}{Theorem} \label{product of dihedral groups}\begin{enumerate}
\item Let $P$ be an arbitrary finite direct product of groups of order 2 and dihedral 2-groups. Then the Dade group of any factor group of $P$ is torsion free.
\item Let $4\leq m_1\leq m_2\leq\cdots\leq m_n$ be a non decreasing sequence of integers. Set $s=\sum_{i=1}^{n-1}\limits(m_i-3)$. Then the torsion part of the Dade group of $\monast{i=1}{n}SD_{2^{m_i}}$ is isomorphic to $(\Z/2\Z)^{2^{s-1}}$ if $m_1<m_n$, and to $(\Z/2\Z)^{2^{s}}$ if $m_1=m_n$.
\item In particular, for any integers $n\geq 1$ and $m\geq 4$, the torsion part of the Dade group of $(SD_{2^m})^{*n}$ is isomorphic to $(\Z/2\Z)^{2^{(n-1)(m-3)}}$.
\end{enumerate}
\end{mth}
\pf First by Example~\ref{D8Q8}
$$D_8\cong\un\oplus 4\cdot \partial C_2$$
in $\mathcal{R}_2$. Now by Corollary~\ref{edge cyclic center}, since the center $Z$ of $D_{16}$ has order 2, and since $D_{16}/Z\cong D_8$
$$D_{16}\cong \un\oplus 4\cdot \partial C_2\oplus \partial D_{16}\mpoint$$
Since for $n\geq 3$, the group $D_{2^n}$ has a center $Z$ of order 2, and since $D_{2^n}/Z\cong D_{2^{n-1}}$, it follows by induction that
\begin{equation}\label{d2n}
D_{2^n}\cong \un\oplus 4\cdot \partial C_2\oplus\mathop{\bigoplus}_{l=4}^n\limits \partial D_{2^l}\mpoint
\end{equation}
Now by Corollary~\ref{edge product}, the product $\partial D_{2^l}\times\partial D_{2^m}$, for $l\leq m$, is isomorphic to $2^{l-3}\cdot\partial D_{2^{m}}$. Moreover $\partial C_2\times \partial P\cong\partial P$ for any 2-group $P$, since $\partial C_2\cong\un$ by Proposition~\ref{edge C2}, it follows that any product of dihedral 2-groups and groups of order~2 is isomorphic to a direct sum of edges of the trivial group, of the edge of the group of order 2, and of edges of dihedral 2-groups.\par
It follows that if $P$ is a direct product of dihedral 2-groups and groups of order 2, then $P$ is isomorphic in $\mathcal{R}_2$ to the direct sum of the trivial group, and some copies of edges of the group of order~2 and the edge of dihedral 2-groups. In other words, if $S$ is a genetic subgroup of $P$, then $N_P(S)/S$ is either trivial, or of order~2, or dihedral. Now the Dade group of dihedral 2-groups is torsion free (by~\cite{cath} Theorem~10.3 (a)), and the Dade groups of the trivial group and the group of order 2 are trivial. It follows that the Dade group of $P$ is torsion free, as well as the Dade group of any quotient of $P$, as is it a direct summand of the Dade group of $P$. This proves Assertion~1.\par
For Assertion~2, set $P=\monast{i=1}{n}SD_{2^{m_i}}$. By Assertion~2 of Corollary~\ref{edge of 2-groups},\vspace{-2ex} 
\begin{equation}\label{pisd}P\cong \prod_{i=1}^m\partial SD_{2^{m_i}}\oplus \prod_{i=1}^nD_{2^{m_i-1}}\mvirg
\end{equation}
in the category $\mathcal{R}_2$. Now an easy induction on $n$, using Corollary~\ref{edge product}, shows that 
$$\prod_{i=1}^n\partial SD_{2^{m_i}}\cong\left\{\begin{array}{cl}
2^{s-1}\cdot\partial C_{2^{m_n-1}}&\hbox{ if }m_1<m_n\\
2^s\cdot \partial SD_{2^{m_n}}&\hbox{ if }m_1=m_n\end{array}\right.$$
where $s=\sum_{i=1}^{n-1}\limits(m_i-3)$.\par
By~\ref{pisd}, this means that in a genetic basis $\mathcal{B}$ of $P$, there are $2^{s-1}$ or $2^s$ subgroups $S$ such that $N_P(S)/S$ is semidihedral, depending on $m_1<m_n$ or $m_1=m_n$, and for the other $S\in\mathcal{B}$, the group $N_P(S)/S$ is trivial, of order~2, or dihedral. \par
The Dade group of a dihedral 2-group is torsion free, and the Dade groups of the trivial group and $C_2$ are trivial. Moreover, the faithful torsion part $\partial D^t(C_{2^m})$ of the Dade group of $C_{2^m}$ is isomorphic to $\Z/2\Z$, if $m\geq 2$ (see \cite{bisetfunctors} Theorem 12.10.3). Similarly, the faithful torsion part $\partial D^{t}(SD_{2^m})$, for $m\geq 4$, is isomorphic to $\Z/2\Z$. This completes the proof of Assertion~2. Assertion~3 is a particular case of Assertion~2.\findemo
\begin{rem}{Remark} \label{units burnside}Let $P$ be a finite product of groups of order 2, and dihedral 2-groups, as in Assertion~1, and let $Q$ be a quotient of $P$. If $T$ is a genetic subgroup of $Q$, then $N_Q(T)/T$ is either trivial, of order 2, or dihedral~: indeed $T$ lifts to a genetic subgroup $S$ of $P$, such that $N_P(S)/S\cong N_Q(T)/T$. It follows in particular that the map
$$\sur{\epsilon}_Q:B^\times(Q)\to\Hom_\Z\big(R_\Q(Q),\F_2\big)$$
introduced in \cite{burnsideunits}, Notation~8.4, is a group isomorphism from the group of units of the Burnside ring of $Q$ to the $\F_2$-dual of $R_\Q(Q)$~: indeed, there are non-negative integers $a$ and $b_i$, for $i\in\{4,\ldots,m\}$ such that
$$Q\cong \un\oplus a\cdot\partial C_2\oplus \dirsum{i=4}^mb_i\cdot \partial D_{2^i}$$
in the category $\mathcal{R}_2$. Then $B^\times(Q)\cong (\F_2)^r$, where $r=1+a+\sum_{i=4}^m\limits b_i$, by \cite{burnsideunits}, Theorem~8.5. Similarly $R_\Q(Q)\cong\Z^r$ (hence $r$ is equal to the number of conjugacy classes of cyclic subgroups of $Q$), so $\Hom_\Z\big(R_\Q(Q),\F_2\big)\cong(\F_2)^r$. As $\sur{\epsilon}$ is injective, it is an isomorphism.
\end{rem}
\begin{mth}{Proposition} Let $m\geq 3$ be an integer. Then for any integer $n$, there is an isomorphism
\begin{equation}\label{D2nk}(D_{2^m})^n\cong \un\oplus (5^n-1)\cdot\partial C_2\oplus\mathop{\bigoplus}_{l=4}^m\limits\frac{(3+2^{l-2})^n-(3+2^{l-3})^n}{2^{l-3}}\cdot\partial D_{2^l}
\end{equation}
in the category $\mathcal{R}_2$.
\end{mth}
\pf Let $\mathcal{S}_p$ denote the full subcategory of $\mathcal{R}_p$ consisting of finite direct sums of edges of Roquette $p$-groups, and let $\Gamma=K_0(\mathcal{S}_p)$ be the Grothendieck group of this category, for relations given by direct sum decomposition. Then Corollary~\ref{simplifiable} shows that $\Gamma$ is a free abelian group, and that two objects of $\mathcal{S}_p$ have the same image in $\Gamma$ if and only if they are isomorphic in~$\mathcal{R}_p$. Moreover, by Corollary~\ref{edge product}, the category $\mathcal{S}_p$ is a tensor subcategory of $\mathcal{R}_p$, and $\Gamma$ is actually a commutative ring.\par
It follows that $\Gamma$ identifies to a subring of the $\Q$-algebra $\Q \Gamma=\Q\otimes_\Z\Gamma$, and that, to prove the proposition, it suffices to check that the two sides of Equation~\ref{D2nk} have the same image in $\Q\Gamma$. Let $c$ denote the image of $\partial C_2$ in~$\Gamma$, and for $l\geq 4$, let $d_l$ denote the image of $\partial D_{2^l}$ in $\Gamma$. By Equation~\ref{d2n}, the image $i_m$ of $D_{2^m}$ in $\Gamma$ is equal to
$$i_m=1+4c+\sum_{l=4}^md_l\mpoint$$
By Corollary~\ref{edge product}, for $4\leq l\leq k$
$$d_l\times d_k = 2^{l-3}d_k\mpoint$$
It follows that the elements $e_l=\dsp\frac{1}{2^{l-3}}d_l$ of $\Q\Gamma$, for $l\geq 4$, are such that
$$\forall l,k,\;4\leq l\leq k,\;\;e_l\times e_k=e_k\mpoint$$
In particular $e_l$ is an idempotent, and the elements
$$f_l=e_l-e_{l+1},\;\hbox{for $4\leq l <m$, and }\;f_m=e_m$$
are orthogonal idempotents of $\Q\Gamma$. With this notation, for $l\geq 4$
$$e_l=f_l+f_{l+1}+\cdots+f_m\mvirg$$
and the element $i_r$ can be written as
\begin{eqnarray*}
i_m&=&1+4c+\sum_{l=4}^md_l\\
&=&1+4c+\sum_{l=4}^m2^{l-3}(f_l+f_{l+1}+\cdots+f_m)\\
&=&1+4c+\sum_{l=4}^m 2(2^{l-3}-1)f_l\mpoint
\end{eqnarray*}
Moreover, it follows from Proposition~\ref{edge C2} that $c\times f_l=f_l$, for $4\leq l\leq m$. Thus
\begin{eqnarray*}
(i_m)^n&=&(1+4c)^n+\sum_{j=1}^n\binom{n}{j}(1+4c)^{n-j}\left(\sum_{l=4}^m 2(2^{l-3}-1)f_l\right)^j\\
&=&(1+4c)^n+\sum_{j=1}^n\binom{n}{j}5^{n-j}\sum_{l=4}^m 2^j(2^{l-3}-1)^jf_l\\
&=&(1+4c)^n+\sum_{l=4}^m\left(\sum_{j=1}^n\binom{n}{j}5^{n-j}2^j(2^{l-3}-1)^j\right)f_l\\
\end{eqnarray*}
\begin{eqnarray*}
(i_m)^n&=&(1+4c)^n+\sum_{l=4}^m\Big(\big(5+2(2^{l-3}-1)\big)^n-5^n\Big)f_l\\
&=&(1+4c)^n+\sum_{l=4}^m\big((3+2^{l-2})^n-5^n\big)f_l\\
&=&(1+4c)^n+\sum_{l=4}^m\big((3+2^{l-2})^n-(3+2^{l-3})^n\big)e_l\mpoint
\end{eqnarray*}
The proposition follows, since $(1+4c)^n=1+(5^n-1)c$, by Equation~\ref{D2n}, and since $e_l=\dsp\frac{1}{2^{l-3}}d_l$.\findemo
\begin{rem}{Remark} The isomorphism~\ref{D2nk} is equivalent to saying that a genetic basis of the group $P=(D_{2^m})^n$ consists of one subgroup $S$ such that $N_P(S)/S\cong \un$ (namely $S=P$), of $5^n-1$ subgroups $S$ such that $N_P(S)/S\cong C_2$, and, for $4\leq l\leq m$, of $\dsp\frac{(3+2^{l-2})^n-(3+2^{l-3})^n}{2^{l-3}}$ subgroups $S$ such that $N_P(S)/S\cong D_{2^l}$.\par
Together with Assertion~3 of Corollary~\ref{edge of 2-groups}, this also gives the structure of genetic bases of the groups $(D_{2^m})^{*n}$, $(SD_{2^m})^{*n}$, $(Q_{2^m})^{*n}$~:
\begin{mth}{Corollary} Let $P$ be one of the groups $D_{2^m}$, $SD_{2^m}$, or $Q_{2^m}$, for $m\geq 4$. Then, for any positive integer $n$, any genetic basis of the group $Q=P^{*n}$ consists~:
\begin{itemize}
\item of one group $S$ such that $N_P(S)/S=\un$ (namely $S=Q$).
\item of $5^n-1$ subgroups $S$ such that $N_P(S)/S\cong C_2$.
\item for $4\leq l\leq m-1$, of $\dsp\frac{(3+2^{l-2})^n-(3+2^{l-3})^n}{2^{l-3}}$ subgroups $S$ such that $N_P(S)/S\cong D_{2^l}$.
\item of $2^{(n-1)(m-3)}$ subgroups $S$ such that $N_P(S)/S$ is isomorphic to $$\left\{\begin{array}{cl}D_{2^m}&\hbox{if $P=D_{2^m}$}\\SD_{2^m}&\hbox{if $P=SD_{2^m}$}\\D_{2^m}&\hbox{if $P=Q_{2^m}$ and $n$ is even}\\Q_{2^m}&\hbox{if $P=Q_{2^m}$ and $n$ is odd}\mpoint\end{array}\right.$$
\end{itemize}
\end{mth}
\end{rem}

\centerline{\rule{5ex}{.1ex}}
\begin{flushleft}
Serge Bouc - CNRS-LAMFA, Universit\'e de Picardie, 33 rue St Leu, 80039, Amiens Cedex 01 - France. \\
{\tt email : serge.bouc@u-picardie.fr}\\
{\tt web~~ : http://www.lamfa.u-picardie.fr/bouc/}
\end{flushleft}
\end{document}